\documentclass[11pt,reqno]{article}
\usepackage{amssymb}
\usepackage{graphicx}
\usepackage{amsmath}
\usepackage[latin1]{inputenc}
\usepackage{amssymb,indentfirst}
\usepackage{times}
\usepackage[T1]{fontenc}
\usepackage{tikz}
\newtheorem{thm}{Theorem}[section]

\newtheorem{lem}[thm]{Lemma}

\def\N{{\mathbb N}}
\def\Z{{\mathbb Z}}
\def\Q{{\mathbb Q}}
\def\R{{\mathbb R}}

\def\bb{\begin}
\def\bc{\begin{center}}
\def\ec{\end{center}}
\def\be{\begin{equation}}
\def\ee{\end{equation}}
\def\ba{\begin{array}}
\def\ea{\end{array}}
\def\bea{\begin{eqnarray}}
\def\eea{\end{eqnarray}}
\def\beaa{\begin{eqnarray*}}
\def\eeaa{\end{eqnarray*}}

\def\hh{\!\!\!\!}

\def\EQ{\hh & = & \hh}
\def\EE{\hh & \equiv & \hh}
\def\LE{\hh & \le & \hh}

\def\AND#1{\hh & #1 & \hh}

\def\al{\alpha}
\def\bt{\beta}
\def\la{\lambda}
\def\La{\Lambda}
\def\vp{\varphi}
\def\ga{\gamma}
\def\Ga{\Gamma}
\def\da{\delta}

\def\e{\varepsilon}
\def\th{\theta}

\def\si{\sigma}
\def\Si{\Sigma}

\def\vpi{\varpi}

\def\bfx{{\bf x}}
\def\bfy{{\bf y}}

\def\bo{{\bf 0}}
\def\bv{{\bf v}}
\def\bfe{{\bf E}}

\def\D{{\mathcal D}}

\def\mcs{{\mathcal S}}

\def\d{\cdot}
\def\dd{\cdots}
\def\oo{\infty}
\def\f{\frac}
\def\pa{\partial}
\def\z{\left}
\def\y{\right}
\def\q{\quad}
\def\qq{\qquad}
\def\uto{\uparrow}

\def\lb{\label}
\def\nn{\nonumber}

\def\rd{\,{\rm d}}
\def\dt{\,{\rm d}t}
\def\ds{\,{\rm d}s}
\def\dr{\,{\rm d}r}
\def\du{\,{\rm d}u}

\def\ck{\check}
\def\tl{\tilde}
\def\tm{\times}
\def\bu{$\bullet$\ }

\def\nm#1{\|#1\|}

\def\x#1{{\rm (\ref{#1})}}

\def\qqf{\qquad \forall}
\def\andq{\quad \mbox{ and } \quad}

\def\ifl{\iffalse}
\def\Proof{\noindent{\bf Proof.} }
\def\qed{\hfill $\Box$ \smallskip}

\def\lan{\la_n}
\def\tlan{\tilde \la_n}
\def\lm{{\ell/m}}
\def\plm{{p,\lm}}
\def\pilm{\pi_\plm}
\def\tlann{\la_{\plm;n}}
\def\tpm{T_p(\mu)}
\def\spm{S_p(\mu)}
\def\upm{U_p(\mu)}

\def\bxy{(\bfx(t),\bfy(t))}
\def\tpt{{2\pi T}}
\def\tps{{2\pi S}}

\renewcommand{\raggedright}{\leftskip=0pt \rightskip=0pt plus 0cm}
\raggedright

\begin{document}

\title
{\bf On the Structure of Periodic Eigenvalues of \\ the Vectorial $p$-Laplacian}

\author{Changjian Liu$^1$, \qq Meirong Zhang$^2$\footnote{Correspondence author.}}

\date{}%

\maketitle

\begin{center}

$^1$ School of Mathematics (Zhuhai), Sun Yat-sen University, \\ Zhuhai, Guangdong 519082, China \\
E-mail: {\tt liuchangj@mail.sysu.edu.cn} \\
\bigskip
$^2$ Department of Mathematical Sciences, Tsinghua University,\\ Beijing 100084, China\\
E-mail: {\tt zhangmr@tsinghua.edu.cn}\\
\end{center}

\begin{abstract}
In this paper we will solve an open problem raised by Man\'asevich and Mawhin twenty years ago on the structure of the periodic eigenvalues of the vectorial $p$-Laplacian. This is an Euler-Lagrangian equation on the plane or in higher dimensional Euclidean spaces. The main result obtained is that for any exponent $p$ other than $2$, the vectorial $p$-Laplacian on the plane will admit infinitely many different sequences of periodic eigenvalues with a given period. These sequences of eigenvalues are constructed using the notion of scaling momenta we will introduce. The whole proof is based on the complete integrability of the equivalent Hamiltonian system, the tricky reduction to $2$-dimensional dynamical systems, and a number-theoretical distinguishing between different sequences of eigenvalues. Some numerical simulations to the new sequences of eigenvalues and eigenfunctions will be given. Several further conjectures towards to the panorama of the spectral sets will be imposed.
\end{abstract}




\tableofcontents

\section{The Problem and the Main Result}
\setcounter{equation}{0} \lb{1st}

This paper is concerned with the structure of periodic eigenvalues of the vectorial $p$-Laplacian. We will reformulate the problem imposed two decades ago by Man\'asevich and Mawhin \cite{MM00} and give our answer in our way.

For simplicity, the period is always taken as $1$. Let $p\in(1,\oo)$ be an arbitrarily given exponent. The vectorial $p$-Laplacian can be introduced as follows. Consider $1$-periodic motions in the $d$-dimensional Euclidean space $\R^d$, $d=1,2,\dd$
    \[
    \bfx =\bfx(t) : \R \to \R^d.
    \]
The kinetic energy and potential energy of the motion are respectively defined by
    \[
    \int_0^1 \|\dot\bfx(t)\|^p \dt, \andq 
    \int_0^1 \|\bfx(t)\|^p \dt. 
    \]
\ifl
    \beaa
    \int_0^1 \|\dot\bfx(t)\|^p \dt \EQ \int_{\R/\Z} \|\dot\bfx(t)\|^p \dt,\\
    \int_0^1 \|\bfx(t)\|^p \dt \EQ \int_{\R/\Z} \|\bfx(t)\|^p \dt.
    \eeaa
\fi
Here $\dot{}=\f{\rd}{\dt}$ and $\|\d\|$ is the Euclidean norm in $\R^d$. Then, under the geometric constraint on $\bfx(t)$
    \[
    \int_0^1 \|\bfx(t)\|^p \dt ={\rm const}\ne 0,
    \]
the Euler-Lagrange Equation for the critical motions $\bfx(t)$ of the kinetic energy is 
    \be\lb{pd}
    (\phi_p(\dot\bfx))\dot{}+ \la \phi_p(\bfx)=\bo.
    \ee
Here $\la\in \R$ is a spectral parameter and the mapping $\phi_p:\R^d\to \R^d$ is as in \x{phip} below. Considering that $\bfx(t)$ is $1$-periodic, i.e. $\bfx(t)$ satisfies the $1$-periodic boundary condition
    \be \lb{bc}
    \z(\bfx(1), \dot\bfx(1)\y)=\z(\bfx(0),\dot\bfx(0)\y),
    \ee
problem \x{pd}-\x{bc} is the eigenvalue problem we will study in this paper. Here, as usual, $\la$ is an eigenvalue of \x{pd}-\x{bc} if Eq. \x{pd} admits some nonzero solutions $\bfx(t)\not \equiv \bo$ satisfying \x{bc}, while solutions $\bfx(t)$ themselves are called eigenfunctions associated with $\la$ or of $\la$. In other words, eigenvalues and eigenfunctions are respectively the critical values and critical points.

The mapping $\phi_p$ in Eq. \x{pd} is
    \be \lb{phip}
    \phi_p:\R^d\to \R^d,\qq \phi_p(\bfx):= \|\bfx\|^{p-2} \bfx.
    \ee
It is an odd self-homeomorphism of $\R^d$ with the inverse
    \(
    \phi_p^{-1} = \phi_q,
    \) \
where $q:=p/(p-1)\in(1,\oo)$ is the conjugate exponent of $p$, and $\phi_p(\bfx)$ has the following homogeneity
    \be \lb{phip1}
    \phi_p(k \bfx) \equiv \phi_p(k) \phi_p(\bfx) \qqf k \in \R, \ \bfx \in \R^d.
    \ee
As an Euler-Lagrange Equation, \x{pd} is a Hamiltonian system and the dimension $d$ is the degree of freedom. When $d=1$ and $d\ge 2$, \x{pd} is simply called the {\it scalar} or the {\it vectorial} $p$-Laplacian in \cite{MM00}, respectively.

The set of all (real) eigenvalues $\la$ of \x{pd}-\x{bc} is denoted by $\Si_{p,d}$, called the spectral set. It is evident that $\Si_{p,d}$ has the trivial eigenvalue $\la_0=0$ with the constant eigenfunctions $\bfe(t;0) \equiv \bv \in \R^d\setminus \{\bo\}$, and other periodic eigenvalues must be positive. Hence we are interested in nontrivial eigenvalues
    \[
    \Si_{p,d}^*:= \Si_{p,d} \setminus \{0\} \subset (0,+\oo),
    \]
and possibly the properties of the associated eigenfunctions $\bfe_p(t)=\bfe_p(t;\la)$. Note that, by embedding $\R^d$ into a higher dimensional space, say into $\R^{d+1}$, one has the trivial inclusions
    \be \lb{sp3}
    \Si_{p,d}^* \subseteq \Si_{p,d+1}^* \qqf d\in \N.
    \ee

For the case $p=2$ (the $2$-Laplacian or the Laplacian) and for arbitrary dimension $d$, by denoting $\bfx=(x_1,\dd,x_d)$, Eq. \x{pd} is
    \[
    \ddot x_i +\la x_i=0, \qq i=1,\dd,d,
    \]
an uncoupled system of linear oscillations. Taking in account of the periodic boundary condition \x{bc}, the spectral set is
    \[
    \Si_{2,d}^* \equiv \{(2n\pi)^2: n \in \N\}.
    \]
Here and henceforth $\N$ 
denotes the set of positive 
integers. 
In particular, $\Si_{2,d}^*$ and $\Si_{2,d}$ are independent of $d$.

However, when the exponent $p\ne 2$, 
it has been found in \cite{del87, MM00} that the periodic spectral sets $\Si_{p,d}^*$ are different for $d=1$ and $d\ge 2$. To describe such a difference and the open problem in \cite{MM00}, let us  first consider the scalar $p$-Laplacian
    \be \lb{p1}
    (\phi_p(\dot x))\dot{}+ \la \phi_p(x)=0.
    \ee
Here $p=2$ is included. Define the number $\pi_p$  by
    \be \lb{pip}
    \pi_p=\frac{2\pi (p-1)^{1/p}}{p\sin (\pi/p)} \in (2,\pi) \q \mbox{for }p\ne 2, \andq \pi_2=\pi.
    \ee
Since Eq. \x{p1} is an integrable equation, it is well-known that
    \be \lb{sp2}
    \Si_{p,1}^* \equiv \{\lan=\la_{n,p}:= (2n\pi_p)^p: n \in \N\}.
    \ee
Moreover, by using the so-called $p$-cosine function $\cos_p(\d)$ in \cite{L95}, the corresponding eigenfunctions of $\lan$ are given by
    \be \lb{sp21}
    E_p(t;\lan) = A \cos_p(2 n\pi_p (t+ \vp_0))\qqf A>0, \ \vp_0\in \R.
    \ee
Results \x{pip}---\x{sp21} for problem \x{p1}-\x{bc} are well recognized in literature.
In particular, we have from \x{sp3} and \x{sp2} that
    \be\lb{sp4}
    \{\lan=(2n\pi_p)^p: n \in \N\}\subset \Si_{p,d}^*\qqf d\ge 2.
    \ee
Moreover, by using \x{sp21}, the eigenvalues $\lan$ of \x{sp4} admit eigenfunctions
    \be \lb{sp2e}
    \bfe_p(t;\lan) \equiv E_p(t;\lan) \bv \qqf \bv \in \R^d, \ \|\bv\|=1.
    \ee
On the other hand, let us just restrict to $d=2$ and define for $n\in \N$,
    \be \lb{sp5} \bb{split}
    \tlan =\tl \la_{n,p} & := (2n\pi)^p, \\
    \tl \bfe_p(t;\tlan) & :=\tl A \bigl(\cos(2n \pi (t+\tl \vp_0)), \sin (2n\pi (t+\tl \vp_0))\bigr) \qqf \tl A>0, \ \tl\vp_0\in \R.
    \end{split}
    \ee
A direct computation shows that $\tlan\in \Si_{p,2}^*$ are also eigenvalues of \x{pd}-\x{bc}, with $\tl \bfe_p(t;\tlan)$ being associated eigenfunctions.

When $p=2$, for any $n\in \N$, one has $\lan\equiv \tlan$ and both $\bfe_2(t;\lan)$ and $\tl\bfe_2(t;\tlan)$ are eigenfunctions. However, whenever $p\ne 2$, though the sets $\{\lan: n\in \N\}$ and $\{\tlan: n\in \N\}$ may intersect, they are different sequences of eigenvalues, cf. \x{pip}. Then, an open problem imposed by Man\'asevich and Mawhin on page 1306 of \cite{MM00} can be stated as
    \[
    {\bf P}_p: \qq \mbox{\it Do we have\q $\Si_{p,2}^* = \{\la_{n,p}: n\in \N\}\cup \{\tl\la_{n,p}: n\in \N\}$?}
    \]
In fact, their problem was originally imposed in arbitrary dimension $d\ge 2$, and more open problems on the further structures of eigenvalues have also been imposed there.

As far as we know, this problem has been left unsolved. Our answer to problem ${\bf P}_p$ is negative. Precisely, we will prove the following theorem.

    \bb{thm} \lb{main}
Let $p\ne 2 $ be arbitrarily given. Then there exists an infinitely many sequence
    \be \lb{vpim}
    \vpi_1, \q \vpi_2, \q \vpi_3, \q \dd \q \vpi_m, \q \dd
    \ee
of positive numbers, depending on $p$, such that
    \be \lb{pimm}
    \f{\vpi_m}{\vpi_{m'}} \not\in \N \qq \mbox{whenever } m\ne m',
    \ee
and, for any $m\in \N$, $\vpi_m$ yields a sequence of $1$-periodic eigenvalues
    \be \lb{evs-s}
    \La_m:=\z\{ \la_{m,n}: = (2n \vpi_m)^p: n\in \N\y\} \subset \Si_{p,2}^*.
    \ee
Hence $\Si_{p,2}^*$ contains infinitely many `different' sequences of eigenvalues as in \x{evs-s}.
    \end{thm}

To make the statements of the theorem be more clear and precise, we give some remarks.

\bu When $m$ is fixed, $\La_m=\{\la_{m,n}\}_{n\in \N}$ is a strictly increasing sequence of eigenvalues tending to $+\oo$ as $n$ goes to $+\oo$.

\bu Result \x{pimm} means that, as subsets of $(0,+\oo)$, whenever $m'\ne m$, the eigenvalue sequence $\La_{m'}$ is neither $\La_m$ itself nor a subsequence of $\La_m$, though $\La_{m'}$ and $\La_m$ may have intersections. By these, it means that when $m$ is different, $\La_m$ is really a different sequence of eigenvalues.


\bu Essentially, the requirements \x{vpim} and \x{pimm} can be stated as that the infinite set $\{\vpi_m: m\in \N\}$ is integer-independent. See Definition \ref{nind}.

\bu We have not appointed any order for the choice of $\vpi_m$'s in \x{vpim}. One can choose, for example, $\vpi_1=\pi_p$ and $\vpi_2=\pi$. These correspond to the two sequences $\{\lan\}$ and $\{\tlan\}$ of eigenvalues observed and stated in the open problem above. 

\bu As a conclusion, Theorem \ref{main} asserts that $\Si_{p,2}^*$ always contain an infinitely many `different' sequence of eigenvalues, each taking the form $(2n \vpi)^p$, $n=1,2,\dd$ In fact, during the whole proofs of Theorem \ref{main}, we will reveal more properties on the spectral sets $\Si_{p,2}^*$, including that for the corresponding eigenfunctions.

Due to the trivial inclusions \x{sp3}, it seems that the bigger $d$ is, the larger the spectral set $\Si_{p,d}$ is. Because of this, let us restrict in this paper to the case $d=2$. Henceforth, $\Si_{p,2}$ and $\Si_{p,2}^*$ are simplified as $\Si_p$ and $\Si_p^*$ respectively.

Let us outline the whole proofs. More or less, the proofs of Theorem \ref{main} are in a constructive way, except at the last step for the verification of integer-independence of the eigenvalues constructed.

In \S \ref{2nd}, we will collect some basic properties on Eq. \x{pd} and its equivalent first-order systems. An important ingredient is that \x{pd} is a completely integrable Hamiltonian system of degree $2$ of freedom, because the energy and the angular momentum are first integrals of \x{pd}. See Lemma \ref{ci}. Hence motions of \x{pd} are described by dynamical systems on $3$-dimensional compact energy levels which can be further foliated using angular momenta. Moreover, by a scaling technique, we will introduce the scaling eigenvalue problems and scaling systems. In particular, the notion of the scaling momenta $\ck \mu \in [0,1]$ for eigenfunctions of problem \x{pd}-\x{bc} will be introduced. See Definitions \ref{sc} and \ref{scss}. For the Laplacian itself, one sees that scaling momenta depend on eigenfunctions, not only on eigenvalues. See Example \ref{se-lan}. All future constructions for eigenvalues will be based on the scaling momenta. For example, in Propositions \ref{prop1}, it will be shown that the two known eigenvalue sequences $\{\lan\}$ and $\{\tlan\}$ will admit some eigenfunctions of the scaling momenta being respectively $0$ and $1$, the possible minimal and maximal scaling momenta.

In \S \ref{3rd}, for the scaling eigenvalue problem on the $3$-dimensional compact energy level, we will deduce a $2$-dimensional reduced dynamical system \x{S2d} on a rectangle $(0,1)\tm (0,\pi)$. In fact, the reduced system has a very simple dynamics, consisting of a single equilibrium and a family of closed orbits surrounding the equilibrium. See Proposition \ref{RDS}. The complete reduction and the complete recovery relations are given in Proposition \ref{RDS2}. For example, the second eigenvalue sequence $\{\tlan\}$ corresponds to the equilibrium of the reduced system. By using the closed orbits of the reduced system, we will introduce two important functions $\tpm$ and $\spm$ using singular integrals. These are used to express the precise dynamical behaviors of the reduced $2$-dimensional and the original $3$-dimensional dynamical systems.

In \S \ref{4th}, we will first use the functions $\tpm$ and $\spm$ to give in Theorem \ref{main1} a complete characterization to the spectral set. Here the momenta $\mu$ ranging from $0$ and $1$ are used. However, as mentioned in the remarks after Theorem \ref{main}, in order to extract really new eigenvalue sequences, we need to verify some integer-independence for the basic eigenvalue sequences like in \x{vpim}. It is apparent that this is a different problem. By the analysis techniques for singular integrals as done in bifurcation theory and limit cycles, we are able to deduce useful information on $\tpm$ and $\spm$ when the scaling momenta $\mu$ are close to $1$. Finally we are able to exploit elementary number theory to construct an infinite integer-independent sequence like in \x{vpim} and complete the proof of Theorem \ref{main}. In particular, as a consequence of the proofs, the sequence
$\{\vpi_m\}_{m\in \N}$ of \x{vpim} for the construction of eigenvalue sequences can be taken so that $\vpi_m\to+\oo$ as $m\to+\oo$. Moreover, each sequence $\La_m$ of eigenvalues will have associated eigenfunctions of scaling momenta tending to the maximal value $1$. For the precise statements of these features, see Corollary \ref{pana1}.

Though the statement of Theorem \ref{main} for the spectral sets $\Si_p^*$, $p\ne 2$, is fairly satisfactory in some sense, the panorama of the complete structures of the spectral sets and sets of scaling momenta is far from being understood in this paper. In order to motivate the future studies to this interesting problem, we will give in \S \ref{sec5} some numerical simulations to our results, but not restricting to large scaling momenta. Both eigenvalues and eigenfunctions will be plotted with some choices of the exponents $p$ and the labelling indexes we are used in Theorem \ref{main1}. Typically, we will arrive at the Lissoajous figures for the periodic motions of eigenfunctions in the configuration space $\R^2$.

Finally, in \S \ref{sec6}, we will first impose some further conjectures on the possible complete structures of the spectral sets $\Si_{p}^*$ and the sets of scaling momenta. Moreover, we will give some explanations to our results from the point of view of the Lagrangian mechanics. It seems to us that the Laplacian here, a special Newtonian equation, is the degenerate case of the $p$-Laplacian as $p\to 2$. In fact, when $p$ is very large or $p$ is very close to $1$, it seems that the structure of eigenvalue will become relatively simpler if the scaling momenta are taken in account. Several references on the eigenvalue problems of different kinds of $p$-Laplacian will be also mentioned at the end of the paper.

\section{Preliminary Results on Eigenvalues and Complete Integrability}
\setcounter{equation}{0} \lb{2nd}

\subsection{Solutions to the vectorial $p$-Laplacian and the equivalent systems}

Let $d=2$ and the exponent $p\in(1,\oo)$ be arbitrarily fixed, including the case $p=2$. Since we need only to study positive eigenvalues of the vectorial $p$-Laplacian \x{pd}-\x{bc}, for convenience, we write $\la$ as $\la^p$ and consider the following equation
    \be\lb{R}
    (\phi_p(\dot\bfx))\dot{}+ \la^p \phi_p(\bfx)=\bo, \qq \bfx \in \R^2.
    \ee
Here $\la>0$ is a positive parameter. By introducing
    \be \lb{y-x}
    \bfy =(y_1,y_2)=\phi_p(\dot\bfx) \in \R^2,
    \ee
which may be called the generalized velocity of the solution $\bfx(t)$, Eq. \x{R} is equivalent to the first-order system
    \be \lb{HS}
    \dot\bfx= \phi_q(\bfy)=\|\bfy\|^{q-2} \bfy,\qq
    \dot\bfy = -\la^p \phi_p(\bfx)= -\la^p \|\bfx\|^{p-2} \bfx,
    \ee
where $(\bfx, \bfy) \in\R^2\tm \R^2 =\R^4$. By introducing $H_\la: \R^4\to \R$ as
    \be \lb{Ham}
    H_\la(\bfx, \bfy):= {\la^p \|\bfx\|^p}/p+ {\|\bfy\|^q}/q,
    \ee
it is well-known that \x{HS} is an autonomous Hamiltonian system of degree $2$ of freedom, because \x{HS} can be written as
    \[
    \dot\bfx= 
    \f{\pa H_\la(\bfx, \bfy)}{\pa \bfy},\qq
    \dot\bfy = 
    -\f{\pa H_\la(\bfx, \bfy)}{\pa \bfx}.
    \]
Hence the Hamiltonian $H_\la(\bfx,\bfy)$ of \x{Ham} is a first integral of system \x{HS}.


Though some of the equations in \x{HS} has singularity, say the second equation of \x{HS} for the case $1<p <2$, solutions of initial value problems to Eq. \x{R} and to system \x{HS} are uniquely and globally defined on $\R$. See \cite{del87} or Lemma 3.1 of \cite{MM00}. Consequently, solutions $\bxy$, $t\in \R$ of \x{HS} are
    \[
    \mbox{either $\z(\bfx(t), \bfy(t)\y)\equiv \bo$ for all $t\in \R$,\qq or $\z(\bfx(t),\bfy(t)\y)\ne \bo$ for all $t\in \R$.}
    \]
Moreover, as $\la>0$ has been assumed, we can conclude from the uniqueness the following equivalence results
    \be \lb{nzs}
    \z(\bfx(t),\bfy(t)\y)\ne \bo\q\forall t\in \R \iff \bfx(t)\not \equiv \bo\iff \bfy(t)\not \equiv \bo.
    \ee
This is a basic observation on solutions which is important in the analysis to eigenvalue problems. See, for example, the proof of Lemma \ref{nz} below.

\subsection{Complete integrability and scaling momenta}

A crucial fact on the Hamiltonian system \x{HS} is its complete integrability.

    \bb{lem} \lb{ci}
System \x{HS} is a completely integrable Hamiltonian system in the sense of Liouville. More precisely, besides the Hamiltonian $H_\la(\bfx, \bfy)$, the (angular) momentum
    \[
    M(\bfx, \bfy): = \bfx \tm \bfy \equiv x_1 y_2 - x_2 y_1 
    \]
is another first integral. Moreover, $H_\la(\bfx, \bfy)$ and $M(\bfx, \bfy)$ are convolutive.
    \end{lem}

\Proof For any solution $\bxy$ of \x{HS}, one has
    \[
    (\bfx \tm \bfy)\dot{} = \dot\bfx \tm \bfy + \bfx \tm \dot\bfy= (\|\bfy\|^{q-2} \bfy) \tm \bfy - \bfx \tm (\la^p \|\bfx\|^{p-2} \bfx) =\bo.
    \]
Hence $M(\bfx,\bfy)$ is also a first integral of \x{HS}.

Moreover, the Poisson bracket of $H_\la$ and $M$  is
    \beaa
    \{H_\la,M\}\EQ\sum_{i=1}^2 \z( \f{\pa H_\la}{\pa x_i} \f{\pa M}{\pa y_i} -\f{\pa H_\la}{\pa y_i} \f{\pa M}{\pa x_i} \y)\\
    \EQ \z( \la^p \|\bfx\|^{p-2} x_1 (-x_2) - \|\bfy \|^{q-2} y_1 y_2 \y) + \z( \la^p \|\bfx\|^{p-2} x_2 x_1 - \|\bfy \|^{q-2} y_2 (-y_1)\y)\\
    \EE 0.
    \eeaa
Hence system \x{HS} is completely integrable in the Liouville sense. The complete integrability has also been used in \cite{HM}. \qed

Associated with any solution $\bxy$ of system \x{HS} are
    \be \lb{HM1}
    h\equiv H_\la\bxy \ge 0,\andq \mu \equiv M\bxy \in\R.
    \ee
They are respectively the $H_\la$-energy and the angular momentum, or simply the energy and momentum of the solution $\bxy$.

From the uniqueness of solutions of initial value problems to the $p$-Laplacian, we have the following observation.

    \bb{lem} \lb{nz}
Let $\bxy$ be any solution of \x{HS} with the zero momentum $\mu=0$. Then there exists some constant unit vector $\bv$ such that
    \be \lb{sol1}
    \bfx(t) \equiv x(t) \bv, \andq \bfy(t) \equiv y(t) \bv,
    \ee
where $(x(t),y(t))$ is a solution to  the following Hamiltonian system of degree $1$ of freedom
    \be \lb{H1}
    \dot x=\phi_q(y),\qq \dot y =-\la^p \phi_p(x).
    \ee
Geometrically, if $\mu=0$, then both $\bfx(t)$ and $\bfy(t)$ are collinear at all times $t\in \R$.
    \end{lem}

\Proof In case $\bfx(t)\equiv \bo$ on $\R$, the results are trivial.

In the following, we assume that $\bfx(t)\not\equiv \bo$ on $\R$. Hence the set
    \[
    D_{\bfx} :=\z\{t\in \R: \bfx(t) \ne \bo \y\}
    \]
is a non-empty open subset of $\R$. Moreover, since $\bfx(t)$ has only non-degenerate zeros, $D_{\bfx}^c$, the compliment of $D_{\bfx}$, is a discrete subset of $\R$.

Let $I$ be any maximal open interval of $D_{\bfx}$. As $\mu=0$, the second equality \x{HM1} is
    \[
    x_1(t) y_2(t) \equiv x_2(t) y_1(t), \qq t\in \R.
    \]
On the interval $I$, this can be solved as
   \be \lb{bfyt}
   \bfy(t) \equiv c(t) \bfx(t) \qqf t\in I,
   \ee
where the scalar function $c(t)$ is
    \[
    c(t)\equiv \f{\bfy(t)\d \bfx(t)}{\bfx(t)\d \bfx(t)}, \qq t\in I.
    \]
Equality \x{bfyt} means that, on $I$, $\bfy(t)$ is parallel to $\bfx(t)$. Then, from the first equation of \x{HS} and equality \x{bfyt}, we have
    \[
    \dot \bfx(t) = \|c(t) \bfx(t)\|^{q-2} c(t) \bfx(t)=: \tl c(t) \bfx(t), \qq t\in I,
    \]
where $\tl c(t) : =  \phi_p(c(t))\|\bfx(t)\|^{q-2}$ is a scalar function on $I$. Hence, with any choice of $t_0\in I$, we obtain
    \be \lb{xti}
    \bfx(t) = \exp\z(\int_{t_0}^t \tl c(s) \ds \y)\bfx(t_0) = x(t) \bv, \qq t\in I,
    \ee
where $\bv=\bfx(t_0)/\|\bfx(t_0)\|$ is a constant unit vector and $x(t)$ is a scalar function on $I$. By using \x{y-x} and result \x{xti} we have obtained, one has also
    \be \lb{yti}
    \bfy(t)=\phi_p(\dot x(t) \bv)= \phi_p(\dot x(t)) \bv=: y(t) \bv, \qq t\in I.
    \ee
Here $y(t)$ is also a scalar function on $I$. As $\bxy$ solves \x{HS}, we have from \x{xti}-\x{yti} that $(x(t),y(t))$ solves \x{H1}, but only restricted to the interval $t\in I$. In particular, $x(t)$ solves the scalar $p$-Laplacian \x{p1}, also only restricted to $t\in I$.

In order to complete the proof, when $D_{\bfx}$ contains more than one maximal interval, we need to show that unit vector $\bv$ in \x{xti}-\x{yti} can be chosen so that it is independent of intervals $I$ of $D_{\bfx}$. To this end, let
    \[
    I_1=(\al,\bt)\andq I_2=(\bt,\ga)
    \]
be two neighboring maximal intervals of $D_{\bfx}$, where
    \(
    -\oo\le \al < \bt < \ga \le +\oo.
    \)
Let us use the common end-point $\bt$ to define the unique solution $x_0(t)$, $t\in \R$ of Eq. \x{p1} satisfying the initial values $(x_0(\bt), \dot x_0(\bt)) =(0,1)$.

Recall that
    \[
    \bfx(t) = x_k(t) \bv_k \q \mbox{ for } t\in I_k, \ k=1,2.
    \]
By the choice of $I_k$, at the common end-point $\bt$ of $I_1$ and $I_2$, one has
    \[
    x_1(\bt-)=0= x_2(\bt+),\q
    \da_1:=\dot x_1(\bt-)\ne 0\andq \da_2:=\dot x_2(\bt+)\ne 0.
    \]
By using the solution $x_0(t)$ defined before, one has from the uniqueness that
    \[
    x_k(t) \equiv \da_k x_0(t), \qq t\in I_k,\ k=1,2.
    \]
Thus, for $t\in I_k$, $k=1,2$, one has
    \[
    \bfx(t) \equiv \da_k x_0(t) \bv_k, \andq
    \dot\bfx(t) \equiv \da_k \dot x_0(t) \bv_k.
    \]
As $\bfx(t)$ and $x_0(t)$ are $C^1$ on $\R$, one has
    \[
    \da_1 \dot x_0(\bt)\bv_1= \dot\bfx(\bt-)= \dot\bfx(\bt+)=\da_2\dot x_0(\bt) \bv_2.
    \]
Since $\dot x_0(\bt)\ne 0$ and $\nm{\bv_1}=\nm{\bv_2}=1$, we obtain
    \[
    \bv_2 = \e \bv_1 \andq \da_2=\e \da_1,
    \]
where $\e=\pm 1$. Hence
    \beaa
    \bfx(t) \EQ \da_1 x_0(t) \bv_1 \q\mbox{ for } t\in I_1,\\
    \bfx(t) \EQ \da_2 x_0(t) \bv_2=\e \da_1 x_0(t) (\e \bv_1) \equiv \da_1 x_0(t) \bv_1\q \mbox{ for } t\in I_2.
    \eeaa
As a whole, by choosing $\bv:= \bv_1$ and $x(t):= \da_1 x_0(t)$, one has
    \[
    \bfx(t) \equiv x(t) \bv \q \mbox{ for } t\in I_1\cup \{\bt\} \cup I_2=(\al,\ga).
    \]
This is the first equality \x{sol1} on a larger interval $(\al, \ga)$ with the same unit vector $\bv$. Inductively, the first equality \x{sol1} can be extended to the whole line $\R$ with the same unit vector $\bv$.

By exploiting the first equality of \x{sol1} on $\R$ we just obtained, we simply set $y(t):= \phi_p(\dot x(t))$ on $\R$. Then
we have the second equality of \x{H1}, and $(x(t),y(t))$ also solves \x{H1} on the whole $\R$. \qed

Due to the homogeneity of system \x{HS}, cf. \x{phip1}, let us introduce the following {\it Scaling System} for \x{HS}
    \be \lb{SS}
    \dot\bfx=\phi_q(\bfy)=\|\bfy\|^{q-2} \bfy,\qq \dot\bfy =-\phi_p(\bfx)= -\|\bfx\|^{p-2} \bfx.
    \ee
It is a special case of \x{HS} with the choice of the parameter $\la=1$. The role of the scaling system \x{SS} is as follows. For any nonzero solution $\bxy$ to system \x{HS}, let $h\equiv H_\la \bxy>0$ and $\mu \equiv M \bxy \in \R$ be its energy and momentum respectively. Let us introduce
    \be \lb{da}
    \da:=\z\{ \ba{ll} +1 & \mbox{ if }\mu \ge 0,\\
    -1 & \mbox{ if }\mu < 0,
    \ea\y.
    \ee
called the {\it Scaling Factor}, and
    \be \lb{SSS}
    \z(\ck \bfx(t),\ck \bfy(t)\y):=\z(\la h^{-1/p}\bfx\z(\da t/\la\y),\da h^{-1/q}\bfy\z(\da t/\la\y)\y),
    \ee
called the {\it Scaling Solution}. Then it is easy to verify that $\z(\ck \bfx(t),\ck \bfy(t)\y)$ is a solution to the scaling system \x{SS}. Moreover, the $H_1$-energy and the momentum of the scaling solution $\z(\ck \bfx(t),\ck \bfy(t)\y)$ are respectively
    \[
    H_1\z(\ck \bfx(t),\ck \bfy(t)\y)\equiv 1, \andq M\z(\ck \bfx(t),\ck \bfy(t)\y)\equiv \ck\mu,
    \]
where
    \be \lb{SM}
    \ck\mu:=|\mu| \la/h.
    \ee
Because of its important role of $\ck \mu$, let us introduce the following definition.

    \bb{defn}\lb{sc}
{\rm The {\it Scaling Momentum} of any nonzero solution $\bxy$ to \x{HS} is defined as $\ck \mu$ in \x{SM}.}
    \end{defn}

\ifl
System \x{SS}, a special case of \x{HS} with the parameter $\la=1$, is called the {\it Scaling System} (of system \x{HS}). The solution $\z(\ck \bfx(t),\ck \bfy(t)\y)$ in \x{SSS} and the constant $\ck\mu$ in \x{SM} are respectively called the {\it Scaling Solution} and {\it Scaling Momentum} (of the solution $\bxy$ to \x{HS}).
\fi

Accordingly, one can say that $\bxy$ has the scaling ($H_1$-)energy $1$. Moreover, the scaling time for the obtention of scaling solutions \x{SSS} can be defined as $s:= \da \la t$. However, we will not use the latter notions too much in the later content.

%
\ifl
    \bb{lem} \lb{scales}
$\bxy$ is a solution of system \x{HS} satisfying $\bxy\in \Ga_{\la,h,\mu}$ iff $\z(\ck \bfx(t),\ck \bfy(t)\y)$ is a solution of the scaling system \x{SS} satisfying $\z(\ck \bfx(t),\ck \bfy(t)\y)\in \Ga_{1,1,\ck \mu}$, where $\ck \mu$ is the scaling momentum of $\bxy$.
    \end{lem}

\Proof We only prove that system \x{HS} can be passed to \x{SS}. For the function in \x{SSS}, one has
    \beaa
    \f{\rd\ck \bfx(t)}{\dt} \EQ h^{-1/p}\da \dot \bfx\z(\da t/\la\y)= h^{-1/p}\da  \phi_q\z(\bfy\z(\da t/\la\y)\y)\\
    \EQ \phi_q\z(\da h^{-1/p(q-1)}\bfy\z(\da t/\la\y)\y)=\phi_q\z(\ck \bfy(t)\y),
    \eeaa

Here we have used \x{phip1} and the equality $p(q-1)=q$. Similarly,
    \beaa
    \f{\rd\ck \bfy(t)}{\dt} \EQ \da h^{-1/q}(\da/\la)\dot\bfy\z(\da t/\la\y)=-h^{-1/q}(1/\la)\la^p \phi_p\z(\bfx\z(\da t/\la\y)\y)\\
    \EQ - \phi_p\z(h^{-1/q(p-1)}\la \bfx\z(\da t/\la\y)\y)=-\phi_p\z( \ck\bfx(t)\y).
    \eeaa
Hence $(\ck \bfx(t),\ck \bfy(t))$ satisfies \x{SS}. Moreover, it is easy to verify that $H_1\z(\ck \bfx(t),\ck \bfy(t)\y)\equiv 1$ and $M\z(\ck \bfx(t),\ck \bfy(t)\y)\equiv \ck\mu$ with $\ck\mu$ being as in \x{SM}.\qed
\fi

One sees from \x{SM} that scaling momenta are non-negative. In fact, one has the following result.

    \bb{lem} \lb{SM0-1}
Any nonzero solution $\bxy$ to \x{HS} has the scaling momentum $\ck \mu\in[0,1]$. Conversely, for any $\ck \mu\in[0,1]$, system \x{HS} has some solution $\bxy$ with the scaling momentum being $\ck \mu$.
    \end{lem}

\Proof For any nonzero solution $\bxy$ to \x{HS}, by the definition \x{HM1} for $h$ and $\mu$, one has from the Cauchy inequality and the H\"older inequality that
    \bea\lb{mu-la}
    |\mu|\EQ |(x_1,x_2)\cdot (y_2,-y_1)|\le \|\bfx\|\cdot \|\bfy\|=\f{1}{\la} \cdot \z(\la\|\bfx\|\cdot \|\bfy\|\y)\nn \\
    \LE \f{1}{\la}\z( \f{\la^p \|\bfx\|^p}{p}+ \f{\|\bfy\|^q}{q}\y)= \f{h}{\la}.
    \eea
Hence \x{SM} implies that $\ck\mu=|\mu| \la/h\in[0,1]$.

In fact, for any $\la>0, \ h>0$ and any $\mu\in \R$, let us define the level set
    \[
    \Ga_{\la,h,\mu}:=\z\{(\bfx,\bfy)\in\R^4: H_\la(\bfx,\bfy)=h, \ M(\bfx,\bfy)=\mu\y\}\ne \emptyset.
    \]
By using the inequality \x{mu-la} and its reversing inequality, one sees that
    \be \lb{Gnp}
    \Ga_{\la,h,\mu}\ne \emptyset\iff \mu\in [-h/\la,h/\la].
    \ee

Finally, for any $\ck \mu\in[0,1]$, one can choose $\mu \in [-h/\la,h/\la]$ such that $|\mu| \la/h=\ck\mu$. Due to the fact \x{Gnp}, let us pick up any point $(\bfx_0,\bfy_0) \in \Ga_{\la,h,\mu}$ as an initial value of system \x{HS} at $t=0$ to obtain a solution $\bxy$ to \x{HS}. From the constructions \x{da}--\x{SM}, one can verify that the scaling momentum of $\bxy$ is just $\ck \mu$ appointed at the beginning. \qed

In fact, when $|\mu|< h/\la$, $\Ga_{\la,h,\mu}$ is diffeomorphic to the $2$-dimensional torus. A famous theorem on completely integrable Hamiltonian systems asserts that all solutions of \x{HS} on $\Ga_{\la,h,\mu}$ are periodic or quasi-periodic.

\ifl
To study solutions of non-zero momenta, we use the scaling technique to simplify the arguments below. To this end, let the parameters be
    \be \lb{para}
    \la>0, \q h>0,\q \mu \in [-h/\la,h/\la], \andq
    \da:=\z\{ \ba{ll} +1 & \mbox{ if }\mu \ge 0,\\
    -1 & \mbox{ if }\mu < 0.
    \ea\y.
    \ee
\fi

    \bb{ex} \lb{se-lan}
{\rm 
(1) Let us consider the first eigenvalue sequence $\{\lan\}$ in \x{sp4} and the associated eigenfunctions $\bfe_p(t;\lan)$ in \x{sp2e}. Since $\bfe_p(t;\lan)$ are along straight lines, the momenta of $\bfe_p(t;\lan)$ are always $0$. By \x{SM}, the scaling momenta are also necessarily $0$. That is, these eigenvalues admit some eigenfunctions of the minimal scaling momentum.

(2) Next we consider the second sequence of eigenvalues and eigenfunctions in \x{sp5}. Transforming to system \x{HS}, for $n\in \N$, one has $\tlan= 2n \pi$ and
    \beaa
    \tl \bfx_n(t)\EQ \tl A \bigl(\cos(2n \pi(t+\tl \vp_0)), \sin (2n\pi(t+\tl \vp_0))\bigr), \\
    \dot{\tl\bfx}_n(t) \EQ 2 n \tl A \bigl(-\sin(2n \pi(t+\tl \vp_0)), \cos (2n\pi(t+\tl \vp_0))\bigr), \\
    \tl\bfy_n(t) \EQ \phi_p(\dot{\tl\bfx}_n(t))\equiv (2 n \tl A)^{p-1} \bigl(-\sin(2n \pi(t+\tl \vp_0)), \cos (2n\pi(t+\tl \vp_0))\bigr).
    \eeaa
Hence the eigenfunction $\tl \bfe_p(t;\tlan) = \z(\tl\bfx_n(t), \, \tl\bfy_n(t)\y)$, as a solution of \x{HS}, has the $\tlan$-energy
    \beaa
    H_{\tlan}\z( \tl\bfx_n(t), \, \tl\bfy_n(t)\y)\EQ \f{1}{p} \tlan^p \|\tl \bfx_n(t)\|^p + \f{1}{q} \|\tl \bfy_n(t)\|^q\\
    \EQ \f{1}{p} (2n\pi)^p \tl A^p + \f{1}{q} (2 n \tl A)^{(p-1)q}\\
    \EQ (2n\pi \tl A)^p=: \tl h_n,
    \eeaa
and the momentum
    \[
    M(\tl\bfx_n(t), \, \tl\bfy_n(t)) = \tl\bfx_n(t)\tm \tl\bfy_n(t)\equiv \tl A (2 n \tl A)^{p-1}= (2n)^{p-1} \tl A^p=:\tl \mu_n.
    \]
By \x{SM}, the scaling momentum is
    \[
\ck \mu_n=|\tl \mu_n| \tlan/\tl h_n = (2n)^{p-1} \tl A^p \d 2n\pi /(2n\pi \tl A)^p \equiv 1.
    \]
That is, these eigenvalues admit some eigenfunctions of the maximal scaling momentum.

(3) Let $p=2$. Then $\pi_2=\pi$ and $\la_{n,2}=\tl\la_{n,2}=2n\pi=:\lan$, $n\in \N$. Due to the linearity, the sums of the corresponding eigenfunctions in \x{sp2e} and in \x{sp5} are still eigenfunctions associated with $\la_n$. More directly, all eigenfunctions can be parameterized as
    \be \lb{E2t}
    \bfe_2(t;\lan)= \bigl(A_1 \sin(\lan(t+\vp_1)), \, A_2 \sin(\lan(t+\vp_2))\bigr),
    \ee
where $A_i\ge 0$, $A_1+A_2>0$, and $\vp_i\in \R$, $i=1,2$. Hence
    \bea \lb{Sme2t}
    h_n \EQ \f{1}{2}\z(\lan^2 \|\bfe_2(t;\lan)\|^2 +\| \dot \bfe_2(t;\lan)\|^2 \y)\equiv \f{1}{2} (A_1^2+A^2_2) \lan^2,\nn\\
    \mu_n \EQ \bfe_2(t;\lan)\tm \dot \bfe_2(t;\lan)\equiv \lan A_1 A_2 \sin \lan(\vp_1-\vp_2), \nn\\
    \ck \mu_n \EQ \f{|\mu_n|\lan}{h_n} \equiv \f{2 A_1 A_2 |\sin \lan(\vp_1-\vp_2)| }{A_1^2+A^2_2}\in[0,1].
    \eea
Result \x{Sme2t} shows that these eigenfunctions in \x{E2t} of the same eigenvalue $\lan$ may have different scaling momenta, which  can actually take all values from $[0,1]$.

(4) Finally, we consider the case $p\ne 2$ such that
    \[
    \f{\pi_p}{\pi} = \f{m_0}{n_0}
    \]
is rational. In this case, the first eigenvalue sequence $\{\lan\}$ in \x{sp4} and the second eigenvalue sequence $\{\tlan\}$ in \x{sp5} will be overlapped, say
    \[
    \la_{n_0}= 2 n_0 \pi_p= 2 m_0 \pi= \tl \la_{m_0}.
    \]
We know from (1) and (2) that such an eigenvalue admits eigenfunctions of the scaling momenta being $0$ and $1$ as well.
\qed}
    \end{ex}

    \bb{rem} \lb{SM-r}
(1) Example \ref{se-lan} shows that scaling momenta of  eigenfunctions of the Laplacian depend on  eigenfunctions themselves, not only on the corresponding eigenvalues.

(2) In Proposition \ref{prop1} below, we will prove that the converse statements  in (1) and (2) of Example \ref{se-lan} are also true in a certain sense.
    \end{rem}

\subsection{A reduction of the eigenvalue problem to the scaling system}

Now we begin to study the spectral set $\Si_p^*$ of problem \x{pd}-\x{bc}. Instead of Eq. \x{pd}, we consider Eq. \x{R} with the spectral parameter changed. Moreover, Eq. \x{R} is equivalent to system \x{HS}, while conditions \x{bc} are transformed to 
    \be \lb{bc1}
    (\bfx(1),\bfy(1))=(\bfx(0),\bfy(0)).
    \ee
For convenience, denote
    \[\bb{split}
    \si_p^*:= & \z\{ \la>0: \mbox{Eq. \x{R} has nonzero solutions $\bfx(t)$ satsifying \x{bc}}\y\}\\
    \equiv & \z\{ \la>0: \mbox{System \x{HS} has nonzero solutions $\bxy$ satisfying \x{bc1}}\y\},
    \end{split}
    \]
cf. the equivalence \x{nzs}. The set $\si_p^*$ is still called the spectral set of the vectorial $p$-Laplacian. It is trivial that $\si_p^*$ differs from $\Si_p^*$ by a power of $1/p$. For example, the eigenvalue sequences \x{sp4} and \x{sp5} are respectively transformed into $\{\lan =2n\pi_p:n\in \N\}$ and $\{\tlan =2n\pi:n\in \N\}$. We will not distinguish such a difference anymore.

Note that the definition \x{SSS} of scaling solutions is essentially irrelevant of the momenta $\mu$. Using the scaling notions, we can characterize $\si_p^*$ in another way.

    \bb{lem} \lb{spic}
For any eigenfunction $\bfe_p(t;\la)$ associated with an eigenvalue $\la\in \si_p^*$ of problem \x{HS}-\x{bc1}, the corresponding scaling solution $\ck \bfe_p(t;\la)$, still simply denoted by $\bxy$, is a solution to the scaling system \x{SS} such that

\bu the scaling energy is
    \be \lb{H11}
    H_1\bxy ={\|\bfx(t)\|^p}/p+ {\|\bfy(t)\|^q}/q \equiv 1,
    \ee

\bu and, the $1$-periodicity of $\bfe_p(t;\la)$ is changed to the $\la$-periodicity of $\bxy$
    \be \lb{P1}
    (\bfx(\la),\bfy(\la)) =(\bfx(0),\bfy(0)).
    \ee

Conversely, solutions $\bxy$ of \x{SS} satisfying \x{H11} and \x{P1} can be transformed into eigenfunctions of problem \x{HS}-\x{bc1} corresponding to the eigenvalue $\la$.
   \end{lem}

\ifl
Suppose that $\bfe_p(t;\la)=\bxy\in \Ga_{\la,h,\mu}$ is any eigenfunction of problem \x{HS}-\x{bc1} associated with $\la$. Note that $(h,\mu)$ depends on $\bfe_p(t;\la)$. As a solution of \x{HS}, $\bfe_p(t;\la)$ has the scaling solution $\ck \bfe_p(t;\la):= \z(\ck \bfx(t), \ck \bfy(t)\y)\in \Ga_{1,1,\ck\mu}$ for some scaling momentum $\ck \mu\in[0,1]$. Moreover,  one sees from \x{SSS} that $\ck\bfe_p(t;\la)$ is a $\la$-periodic solution of the scaling system \x{SS}.
\fi

Due to these facts, for simplicity, let us introduce the following notions.

    \bb{defn} \lb{scss}
{\rm For an eigenfunction $\bfe_p(t;\la)$ of problem \x{HS}-\x{bc1} associated with an eigenvalue $\la\in \si_p^*$,  $\bxy:= \ck \bfe_p(t;\la)$ and $\mu:= M\bxy \in[0,1]$ are respectively called the {\it scaling eigenfunction} and the {\it scaling momentum} as before. Moreover, the number $\la$ is now called the {\it scaling period} (of $\bfe_p(t;\la)$).}
    \end{defn}

For example, in Example \ref{se-lan}, eigenfunctions $\bfe_p(t;\lan)$ have the scaling eigenfunctions
    $$
    \ck\bfe_p(t;\lan)=\bigl(\bv \cos_p(t+\vp_0),\, \bv \phi_p\z(\cos'_p(t+\vp_0) \y)\bigr),
    $$
the scaling momenta $0$, and the scaling periods $\lan=2n\pi_p$. Eigenfunctions $\tl\bfe_p(t;\lan)$ have the scaling eigenfunctions
    \[
    \ck{\tl\bfe}_p(t;\lan)=\bigl(\z(\cos(t+\vp_0),\, \sin(t+\vp_0)\y), \ \z(-\sin(t+\vp_0),\, \cos(t+\vp_0)\y)\bigr),
    \]
the scaling momenta $1$, and the scaling periods $\tlan=2n\pi$.

\ifl

    \bb{lem}
$\la\in \si_p^*$ iff the scaling system \x{SS} admits some solution $\bxy\in \Ga_{1,1,\mu}$ for some $\mu\in[0,1]$ such that it is $\la$-periodic, i.e., iff the following problem \x{S00}-\x{P1} has some solution $\bxy$
    \be \lb{S00}
    \z\{ \ba{l}
    \dot\bfx(t)= \|\bfy(t)\|^{q-2} \bfy(t),\q \dot\bfy(t) = -\|\bfx(t)\|^{p-2} \bfx(t),\\
    {\|\bfx(t)\|^p}/p+ {\|\bfy(t)\|^q}/q\equiv 1, \\
    (x_1(t),x_2(t))\d (y_2(t),-y_1(t))\equiv \mu,
    \ea\y.
    \ee
   \end{lem}

Here, as before, $\bxy$, $\mu$, and $\la$ are respectively the scaling eigenfunction, momentum, and period.
The following analysis is based on the characterization \x{S00}-\x{P1} for $\si_p^*$.\fi

Basing on Lemma \ref{spic}, we can now give a complete explanation to the eigenvalue sequences in \x{sp4} and \x{sp5}.

    \bb{prop} \lb{prop1}
(1) Let $\la\in \si_p^*$. Then $\la$ is equal to $\lan$ for some $n\in \N$ iff there are some eigenfunctions associated with $\la$ of the scaling momentum $0$.

(2) Let $\la\in \si_p^*$. Then $\la$ is equal to $\tlan$ for some $n\in \N$ iff there are some eigenfunctions associated with $\la$ of the scaling momentum $1$.
    \end{prop}

\Proof The necessity parts of the proposition have been explained in Example \ref{se-lan}.

For the sufficiency part of (1), suppose that $\la\in\si_p^*$ admits some eigenfunctions of the scaling momentum $0$. We can now work  on solutions of problem \x{SS}-\x{H11}-\x{P1}. Since $\mu=M\bxy =0$, it follows from Lemma \ref{nz} that $\bxy \equiv (x(t) \bv, y(t) \bv)$ for some unit vector $\bv$. Thus \x{SS}, \x{H11} are \x{P1} respectively 
    \[\bb{split}
    &\dot x=\phi_q(y),\qq \dot y =- \phi_p(x),\\
    & |x(t)|^p/p + |y(t)|^q/q 
    \equiv 1,\\
    & (x(t+\la),y(t+\la)) \equiv (x(t),y(t)).
    \end{split}
    \]
By solving this scalar system, we have already known that the scaling period $\la$ must be $\lan=2n\pi_p$ for some $n\in \N$.

For the sufficiency part of (2),  the main idea of the proof is similar. Let us just consider solutions $\bxy$ to problem \x{SS}-\x{H11}-\x{P1} with
    \[
    M\bxy =(x_1(t),x_2(t))\d (y_2(t),-y_1(t))\equiv 1.
    \]
Instead of Lemma \ref{nz}, we can go back to the proof of \x{mu-la}. Note from \x{H11} and \x{P1} that
    \[
    (y_2(t),-y_1(t))\equiv (x_1(t),x_2(t)), \andq \|\bfx(t)\|=\|\bfy(t)\| \equiv 1.
    \]
Then system \x{SS} for $\bxy$ is simply reduced to the linear system
    \[
    \dot\bfx=  \bfy,\qq \dot\bfy = -\bfx.
    \]
By solving this system, the scaling period $\la$ must be $\tlan=2n\pi$ for some $n\in \N$.\qed


\section{Dynamics of the Reduced Planar Dynamical System}
\setcounter{equation}{0} \lb{3rd}

\subsection{The reduced planar dynamical system}

In Lemma \ref{spic}, we have using the scaling problem \x{SS}-\x{H11}-\x{P1} to characterize the spectral set $\si_p^*$. Note that \x{SS}-\x{H11} is a $3$-dimensional dynamical system on the energy level $H_1^{-1}(1)\subset \R^4$. In this section we will further reduce this to a planar dynamical system which is also integrable. This will be used in the next section, by combining with the boundary condition \x{P1}, to characterize all eigenvalues $\la\in \si_p^*$.

We will prove in Proposition \ref{RDS2} below that \x{SS}-\x{H11}, considered as a dynamical system, can be completely reduced to a planar dynamical system. To describe the reduced system, let us first introduce a planar vertical strip
    \[
    \D:=(0,1)\tm \R.
    \]
Define the functions $F(r)=F_p(r):(0,1) \to \R$ and $G(r)= G_p(r):(0,1) \to \R$ by
    \bea\lb{Fr}
    F(r)\AND{:=} (p r)^{1/q} (q (1-r) )^{1/p}, \\ 
    \lb{Gr}
    G(r) \AND{:=} (p r)^{1/q} /(q (1-r))^{1/q}-(q (1-r))^{1/p}/(p r)^{1/p} \nn \\
    \EE p^{1/q} q^{1/p}\f{r-1/p}{r^{1/p}(1-r)^{1/q}}.
    \eea
It is elementary to verify the following properties for these functions:

\bu $F(r), \ G(r)\in C^\oo(0,1)$;

\bu $G(r)$ is strictly increasing in $(0,1)$ and $G(0+)=-\oo$, $G(1/p)=0$ and $G(1-) = +\oo$; and

\bu $F(r)$ is strictly positive on $(0,1)$, and $F(0+)=F(1-)=0$. See Figure \ref{qprf}.

By using these, let us introduce the following $2$-dimensional dynamical system
    \be \lb{S2d}
    \dot r= F(r) \cos \th,\q
    \dot \th = G(r) \sin \th,\qq (r,\th) \in \D.
    \ee
We give some properties on the vector field of \x{S2d} and its dynamics on $\D$, most of which are evident.

\bu The vector field $(F(r)\cos \th, G(r)\sin \th)$ of \x{S2d} is smooth in the strip $\D$ and $2\pi$-periodic in $\th$.

\bu For any $k\in \Z$, the horizontal segment $H_k:=\{(r,\th): r\in(0,1), \ \th= k\pi\}= (0,1)\tm \{k\pi\}$ is an invariant set of \x{S2d}. Moreover, on $H_k$, system \x{S2d} is determined by the following ODE
    \be \lb{eq-vk}
    \dot r= (-1)^k F(r).
    \ee

\bu System \x{S2d} is invariant under the transformation
    \be \lb{Mat}
    {\mathcal T}:\q t\to -t, \q r\to r, \q \th \to \th+\pi.
    \ee

\bu As $G(0+)=-\oo$ and $G(1-) = +\oo$, the vertical lines $V_0:= \{0\}\tm \R$ and $V_1:=\{1\}\tm \R$ can be understood as singularities of system \x{S2d}.

Due to these properties, especially the invariance results \x{eq-vk} and \x{Mat}, we need only to consider system \x{S2d} in the following fundamental domain
    \[
    \D_0:=\z\{(r,\th): r\in (0, 1), \ \th\in [0,\pi)\y\}=(0,1)\tm [0,\pi),
    \]
which is an invariant rectangle of \x{S2d}. To describe the dynamics of \x{S2d} on $\D_0$, let us introduce the third function $Q(r)=Q_p(r):(0,1) \to \R$ by
    \be \lb{Qr}
    Q(r):=(p r)^{1/p} (q (1-r) )^{1/q}. 
    \ee
It is easy to verify that

\bu $Q(r)$ is a smooth, positive, strictly convex function of $r\in(0,1)$; and

\bu $Q(0+)=Q(1-)=0$, and $\max_{r\in(0,1)} Q(r)=Q(1/p)= 1$.

For example, when $p=q=2$,
    \be \lb{fp1}
    F_2(r)=Q_2(r) \equiv 2 \sqrt{r(1-r)}, \andq G_2(r) \equiv \f{2r-1}{ \sqrt{r(1-r)}}.
    \ee
In the later content, we are mainly using the functions $F_p(r)$ and $Q_p(r)$. For their graphs, see Figure \ref{qprf}.

\begin{figure}[httb]
\bc\resizebox{8cm}{6cm}{\includegraphics{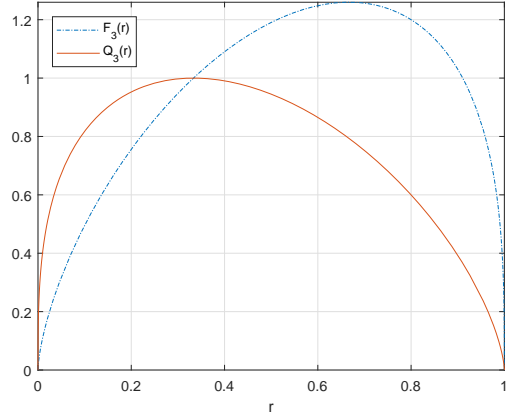}}\ec
\caption{The functions $F_p(r)$ and $Q_p(r)$, where $p=3$. } \lb{qprf}
\end{figure}

From the defining equalities \x{Fr}, \x{Gr} and \x{Qr}, it is easy to establish the following connection
    \[
    Q'(r) F(r) + Q(r) G(r)\equiv 0\qqf r\in(0,1).
    \]
Then, for any solution $(r(t),\th(t))\in \D_0$ of \x{S2d}, we have
    \beaa
    \f{\rd}{\dt} \z(Q(r)\sin\th\y)\EQ Q'(r) \dot r \sin \th + Q(r) (\cos \th) \dot \th\\
    \EQ Q'(r) F(r) \cos \th \sin \th + Q(r) G(r) \sin \th \cos \th \\
    \EE 0.
    \eeaa
Hence the function $Q(r)\sin\th: \D_0 \to \R$ is a first integral of \x{S2d}, which has the range $[0,1]$. Thus there exists $\mu\in[0,1]$ such that
    \be \lb{Qth}
    Q(r(t))\sin\th(t) \equiv \mu.
    \ee

At this moment, we simply understand $\mu\in[0,1]$ as a parameter and use \x{Qth} to define the set
    \be \lb{Cmu}
    C_\mu: =\z\{ (r,\th)\in \D_0: Q(r)\sin\th =\mu\y\}.
    \ee
Then $C_0 =(0,1)\tm \{0\}$ is the bottom side of $\D_0$, and $C_1=\{(1/p,\pi/2)\}$ is a single point. For $\mu\in(0,1)$, $C_\mu$ is a closed curve inside $\D_0$. In fact, let $r_\pm =r_\pm(\mu) =r_{\pm,p}(\mu)\in (0,1)$ be the two solutions of the equation
    \be \lb{Gmu}
    Q_p(r) = \mu, \qq 0 < r_- < 1/p < r_+ < 1.
    \ee
See Figure \ref{qmurf}. Then $C_\mu$ can be parameterized as two curves
    \be \lb{ths}
    C_\mu^-: \q \th =\arcsin \z(\mu/Q(r)\y), \andq C_\mu^+: \q \th =\pi-\arcsin \z(\mu/Q(r)\y),
    \ee
with parameter $r\in[r_-,r_+]$. Hence $C_\mu$ is a simple closed curve. Moreover, $\D_0$ is foliated  as
    \[
    \bigcup_{\mu\in [0,1]} C_\mu=\D_0.
    \]

\begin{figure}[httb]
\bc\resizebox{8cm}{6cm}{\includegraphics{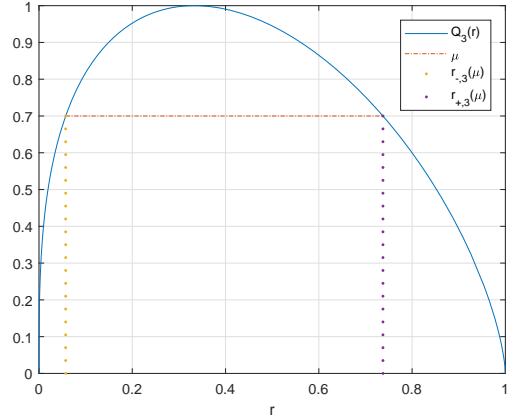}}\ec
\caption{Constructions of $r_{\pm,p}(\mu)$, where $p=3$. } \lb{qmurf}
\end{figure}

Due to \x{Qth}, each $C_\mu$ is actually an orbit of system \x{S2d}. The complete dynamics of the reduced dynamical system \x{S2d} on $\D_0$ is stated in Proposition \ref{RDS}. For the phase portrait, see Figure \ref{fig1}.

\begin{figure}[httb]
\bc\resizebox{8cm}{6cm}{\includegraphics{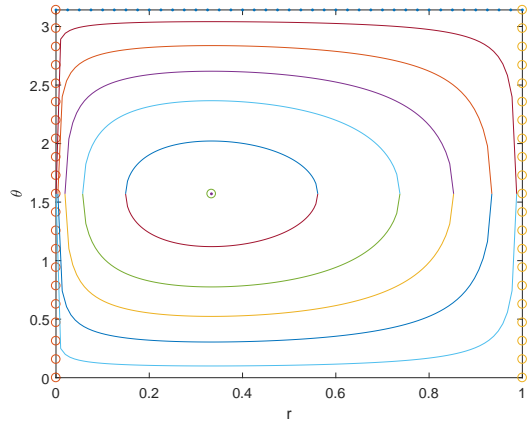}}\ec
\caption{The phase portrait of the reduced system \x{S2d}, where $p=3$. }
\lb{fig1}
\end{figure}

    \bb{prop} \lb{RDS}
The dynamics of system \x{S2d} on $\D_0$ is as follows.

(1) System \x{S2d} is integrable on $\D_0$ in sense that the function $Q(r)\sin\th$ is a first integral.

(2) More precisely;

\bu $C_1$ is the unique equilibrium of \x{S2d},

\bu for any $\mu\in(0,1)$, $C_\mu$ is a nonconstant periodic orbit of \x{S2d}, along which  solutions of \x{S2d} are going anti-clockwise; and

\bu on the bottom segment $C_0$, solutions of \x{S2d} are going right-forward.
    \end{prop}

Now let us establish the relation between the dynamical behaviours of system \x{SS}-\x{H11} and  the reduced system \x{S2d}.

Let $\bxy$ be any solution of \x{SS}. Then $H_1\bxy \equiv 1$, and
    \be \lb{nzs1}
    M\bxy \equiv \mu \in[0,1].
    \ee

To motivate the reduction, let us first consider in \x{nzs1} the case $\mu \in(0,1]$. In this case, results \x{nzs} can be improved as
    \[
    \|\bfx(t)\|>0, \andq \|\bfy(t)\|>0 \qqf t\in\R.
    \]
As a result, $\bfx(t)$ and $\bfy(t)$ can be written as in the following polar coordinates
    \be \lb{PC}
    \bfx = x_1+ i  x_2 = (p r)^{1/p} e^{i \vp}\andq \bfy= y_1+ i  y_2 =(q s)^{1/q} e^{i \psi},
    \ee
where $r>0$ and $s>0$ are respectively the potential and the kinetic energies of the motions $\bxy$. More restrictions on $(r,s, \vp,\psi)$ will be deduced later. By \x{PC}, we have
    \beaa
    \dot\bfx \EQ
    \z( (p r)^{-1/q}\dot r +i (p r)^{1/p} \dot \vp\y) e^{i \vp}, \\
    \dot\bfy \EQ
    \z( (q s)^{-1/p} \dot s +i (q s)^{1/q} \dot \psi \y)e^{i \psi}, \\
    \phi_q(\bfy) \EQ (q s)^{1/p} e^{i \psi},\\
    \phi_p(\bfx) \EQ (p r)^{1/q} e^{i \vp}.
    \eeaa
Thus system \x{SS} becomes
    \beaa
    \z\{\ba{rcl}(p r)^{-1/q}\dot r +i (p r)^{1/p} \dot \vp \EQ (q s)^{1/p} e^{i (\psi-\vp)},\\
    (q s)^{-1/p} \dot s +i (q s)^{1/q} \dot \psi \EQ -(p r)^{1/q} e^{i (\vp-\psi)}.
    \ea\y.
    \eeaa
Taking the real and the imaginary parts, we arrive at
    \be \lb{S1}
    \z\{\ba{rcl} \dot r \EQ (p r)^{1/q} (q s)^{1/p} \cos (\psi-\vp),\\
    \dot s \EQ -(p r)^{1/q} (q s)^{1/p} \cos (\psi-\vp) ,\\
    \dot \vp \EQ (p r)^{-1/p} (q s)^{1/p} \sin (\psi-\vp),\\
    \dot \psi \EQ (p r)^{1/q} (q s)^{-1/q} \sin (\psi-\vp).
    \ea\y.
    \ee
Thus let us introduce
    \be \lb{th-def}
    \th:=\psi-\vp\in \R,
    \ee
i.e., $\th(t)$ is the phase difference between the motion $\bfx(t)$ and its (generalized) velocity $\bfy(t)=\phi_p(\dot \bfx(t))$.

Due to $H_1\bxy \equiv 1$, in the coordinates \x{PC}, one has
    \be \lb{rr11}
    r\in(0,1) \andq s\equiv 1-r\in(0,1).
    \ee
By equality \x{rr11} and the definition \x{th-def} of $\th$, the first equation of \x{S1} has been stated as the first equation of the reduced system \x{S2d}, where the resulted function $F(r)$ is as in \x{Fr}. Moreover, from  the latter two equations of \x{S1} we have
    \[
    \f{\rd\th}{\dt} = \dot \psi-\dot \vp=\z((p r)^{1/q} (q s)^{-1/q} -(p r)^{-1/p} (q s)^{1/p}\y)\sin (\psi-\vp)\equiv G(r) \sin \th,
    \]
by using \x{rr11} and the function $G(r)$ in \x{Gr}. This gives the second equation of \x{S2d}. Hence we have proved that $(r,\th)$ satisfies the reduced system \x{S2d}.

To see that $\th(t)$ can be chosen to be in the interval $(0,\pi)$, let us notice from \x{nzs1} and \x{PC} that
    \bea\lb{mu2}
    \mu \EE x_1(t) y_2(t) - x_2(t) y_1(t) \nn\\
    \EQ \z((p r)^{1/p}\cos\vp\y) \z((q (1-r))^{1/q} \sin \psi\y)- \z((p r)^{1/p}\sin\vp\y) \z((q (1-r))^{1/q} \cos \psi\y) \nn\\
    \EQ (p r)^{1/p}(q (1-r))^{1/q} \sin(\psi-\vp) \nn\\
    \EE Q(r)\sin\th.
    \eea
As we are considering $\mu\in (0,1]$, \x{mu2} shows that $\sin \th(t)\ne 0$ and therefore there exists some $k\in \Z$ such that $\th(t) \in (k\pi,(k+1)\pi)$. Hence, by applying ${\mathcal T}$ in \x{Mat}, one can shift $\vp(t)$ and $\psi(t)$ so that $\th(t) \in (0,\pi)$. Moreover, \x{mu2} also shows that the parameter $\mu$ in \x{Qth} is just the (scaling) momentum of the solution $\bxy$ of system \x{SS}-\x{H11}  in the present case. Hence $\bxy$ can correspond to the orbit $C_\mu$ of the reduced system \x{S2d}.

For the remaining case $\mu=0$, we have obtained \x{sol1}. See Lemma \ref{nz}. In the coordinates \x{PC}, there exists $\vp_0\in\R$ such that
    \[
    \bfx(t) \equiv x(t) e^{i \vp_0}, \andq \bfy(t) \equiv y(t) e^{i \vp_0}.
    \]
These solutions can correspond to the orbit $C_0$ of the reduced system \x{S2d}.

Conversely, let $(r(t), \th(t))\in C_\mu \subset \D_0$ be any solution of the reduced system \x{S2d} for some $\mu\in[0,1]$. Let us define $s(t)$, $\vp(t)$ and $\psi(t)$ using
    \be \lb{vpeq} \bb{split}
    s(t) & := 1-r(t),\\
    \f{\rd\vp(t)}{\dt} & := (p r(t))^{-1/p} (q (1-r(t)))^{1/p} \sin \th(t),\\
    \psi(t)& := \th(t) +\vp(t).
    \end{split}
    \ee
It can be verified that $(r(t),s(t),\vp(t),\psi(t))$ solves system \x{S1}. Then, via \x{PC}, this gives a solution $\bxy$ of system \x{SS}. Moreover, we have already seen that $\bxy$ must have the momentum $\mu$.

These results can be stated as follows.

    \bb{prop}\lb{RDS2}
System \x{SS}-\x{H11} is equivalent to the reduced system \x{S2d} on $\D_0$. More precisely, for any $\mu\in[0,1]$, solutions $\bxy$ of \x{SS}-\x{H11} with the scaling momentum $\mu$ correspond to the orbit $C_\mu$ of \x{S2d}.
    \end{prop}

    \bb{rem}\lb{mu01}
(1) Using the reduced dynamics, Proposition \ref{prop1} asserts that eigenvalues $\{\lan\}$ and $\{\tlan\}$ admit eigenfunctions corresponding to the bottom side $C_0$ and to the equilibrium $C_1$ respectively.

(2) Due to Proposition \ref{RDS2}, when system \x{S2d} is used to obtain new eigenvalues of the $p$-Laplacian different from $\{\lan\}$ and $\{\tlan\}$, it is enough to consider momenta $\mu\in(0,1)$. In this sense, we can consider the reduced system \x{S2d} in the following pinched domain
    \[
    \D_0^*:= \z\{ (r,\th): r\in(0,1), \ \th\in (0,\pi)\y\}\setminus \{(1/p,\pi/2)\} \equiv \bigcup_{\mu\in (0,1)} C_\mu,
    \]
whose dynamics is simply a family of nonconstant periodic orbits.
    \end{rem}

\subsection{Construction of the functions $\tpm$ and $\spm$}

Because of Proposition \ref{prop1} and Remark \ref{mu01}, in order to use the scaling system \x{SS}-\x{H11} and the reduced system \x{S2d} to find new eigenvalues, we need only to consider the scaling momenta $\mu\in(0,1)$.

For any $\mu\in(0,1)$, a solution $\bxy$ of problem \x{SS}-\x{H11} of the scaling momentum $\mu$ corresponds to nonconstant periodic orbit $C_\mu=\{(r(t),\th(t))\}\subset \D_0^*$ of system \x{S2d}.

\bu Let us use $\tpt=\tpt_{p}(\mu)>0$ to denote the minimal period of $(r(t),\th(t))$. Then
    \be \lb{Tdef}
    (r(t+\tpt),\th(t+\tpt)) \equiv (r(t),\th(t)).
    \ee
Here, for later convenience and normalization, a factor $2\pi$ is added to the minimal period.

\bu With such a $\tpt$-periodic solution $(r(t),\th(t))$, Eq. \x{vpeq} shows that $\dot\vp(t)$ is a positive $\tpt$-periodic function. Let us define
    \be \lb{spmu}\bb{split}
    S =& S_{p}(\mu):= \f{1}{2\pi}(\vp(\tpt)-\vp(0))\\
    \equiv & \f{1}{2\pi}\int_0^{\tpt} (p r(t))^{-1/p} (q (1-r(t)))^{1/p} \sin \th(t) \dt>0.
    \end{split}
    \ee
Since $\psi(t)=\vp(t)+\th(t)$, where $\th(t)$ is $\tpt$-periodic, one also has
    \be \lb{Sdef}
    S = 
    \f{1}{2\pi}(\psi(\tpt)-\psi(0)).
    \ee
The role of $S$ is that $\vp(t)$ and $\psi(t)$ satisfy
    \be \lb{vpt}
    \vp(t+\tpt) \equiv \vp(t)+ \tps,\andq \psi(t+\tpt) \equiv \psi(t)+ \tps.
    \ee

Let can express $T=\tpm$ and $S=\spm$ as singular integrals. Notice from \x{ths} that $C_\mu^-$ and $C_\mu^+$ are symmetric with respect to $\th=\pi/2$. By a translation of time, one can assume that
    \[
    (r(0),\th(0))=(r_-,\pi/2)\in C_\mu^-.
    \]
Then $r(t)$ and $\th(t)$ satisfy
    \[
    r(\tpt-t)\equiv r(t)\andq \th(\tpt-t) \equiv \pi-\th(t).
    \]
Moreover, one has $r(0)= \min_t r(t) = r_-$,
    \[
    r(\pi T)=\max_t r(t) =r_+,\andq \th(\pi T)=\pi/2.
    \]
See the phase portrait in Figure \ref{fig1}. For $t\in[0,\pi T]$, one has

\bu $\th(t) \in(0,\pi/2]$, $(r(t), \th(t))\in C_\mu^-$,

\bu it follows from \x{Cmu} that
    \be \lb{tht1}
    \sin \th(t) = \f{\mu}{Q(r(t))},
    \ee

\bu with \x{tht1}, the first equation of system \x{S2d} gives
    \be \lb{E5}
    \f{\dr}{\dt}=\dot r=  F(r)\sqrt{1 -\z(\f{\mu}{Q(r)}\y)^2}= \f{F(r)}{Q(r)}\sqrt{Q^2(r) -\mu^2}.
    \ee

\noindent Integrating equation \x{E5}, we obtain the minimal period of closed orbit $C_\mu$. By considering the factor $2\pi$, $T$ is as in \x{Tmu}.

    \bb{lem}\lb{Tmu0}
For $\mu\in(0,1)$, $T=T(\mu)=\tpm$ is given by
    \bea \lb{Tmu}
    T(\mu)\EQ \f{1}{\pi} \int_{r_-(\mu)}^{r_+(\mu)} \f{Q(r)}{F(r)\sqrt{Q^2(r) -\mu^2}}\dr.
    \eea
    \end{lem}

Next let us work out the expression of $S$. For $t\in (0,\pi T]$, we have from \x{vpeq}, \x{tht1} and \x{E5} that
    \beaa
    \rd\vp\EQ (p r)^{-1/p} (q (1-r))^{1/p} \sin \th\dt\\
    \EQ (p r)^{-1/p} (q (1-r))^{1/p}\f{\mu}{Q(r)}\dt\\
    \EQ (p r)^{-1/p} (q (1-r))^{1/p}\f{\mu}{Q(r)} \f{Q(r)}{F(r)\sqrt{Q^2(r) -\mu^2}}\dr\\
    \EE \f{\mu \dr}{p r \sqrt{Q^2(r) -\mu^2}}.
    \eeaa
By considering the factor $2\pi$, we obtain
    \be \lb{Smu0}
    S=\f{\mu}{\pi} \int_{r_-(\mu)}^{r_+(\mu)} \f{\dr}{p r \sqrt{Q^2(r) -\mu^2}}.
    \ee
If the equation for $\psi(t)$ is used, one can obtain another expression
    \[
    S =\f{\mu}{\pi} \int_{r_-(\mu)}^{r_+(\mu)} \f{\dr}{q(1- r) \sqrt{Q^2(r) -\mu^2}}.
    \]
Then, by taking a weighted average, we can obtain an expression for $S$ in a symmetric form.

    \bb{lem} \lb{Smu-e}
For $\mu\in(0,1)$, $S=S(\mu)=S_p(\mu)$ is given by
    \be \lb{Smu}
    S(\mu)= \f{\mu}{\pi}\int_{r_-(\mu)}^{r_+(\mu)} \f{\dr}{pr \d q(1-r) \sqrt{Q^2(r) -\mu^2}}.
    \ee
    \end{lem}

    \bb{ex} \lb{TSp2}
{\rm \bu Let us consider the case $p=q=2$. One has
    \[
    F_2(r)= Q_2(r)\equiv \sqrt{ 4r-4r^2}.
    \]
See \x{fp1}. Hence the solutions of the corresponding equation \x{Gmu} are
    \[
    r_{\pm,2}(\mu) =\f{1}{2} \z(1\pm \sqrt{1-\mu^2}\y).
    \]
As $\f{Q_2(r)}{F_2(r)}\equiv 1$, we have from \x{Tmu} and \x{Smu0} that
    \bea \lb{T2mu}
    T_2(\mu) \EQ \f{1}{\pi} \int_{r_{-,2}(\mu)}^{r_{+,2}(\mu)} \f{\dr}{\sqrt{4 r-4 r^2-\mu^2}}\equiv \f{1}{2} \qqf \mu \in(0,1),\\
    \lb{S2mu}
    S_2(\mu) \EQ \f{\mu}{\pi} \int_{r_{-,2}(\mu)}^{r_{+,2}(\mu)} \f{\dr}{2r\sqrt{4 r-4 r^2-\mu^2}}\equiv \f{1}{2}\qqf \mu \in(0,1).
    \eea

\bu Note that \x{Tmu} and \x{Smu} are singular integrals. For example, for $p=3$, $r=r_{\pm,3}(\mu)$ are determined by
    \[
    Q_3(r) = \f{3 r^{1/3}(1-r)^{2/3}}{2^{2/3}}=\mu,
    \]
which is essentially a cubic equation, while
    \beaa
    T_3(\mu) \EQ \f{1}{\pi}\int_{r_{-,3}(\mu)}^{r_{+,3}(\mu)}\f{ (1-r)^{1/3}}{(2 r)^{1/3}\sqrt{\f{9}{2^{4/3}} r^{2/3} (1-r)^{4/3}-\mu^2} } \rd r,\\
    S_3(\mu) \EQ \f{1}{\pi}\int_{r_{-,3}(\mu)}^{r_{+,3}(\mu)}\f{\mu}{3 r\sqrt{\f{9}{2^{4/3}} r^{2/3} (1-r)^{4/3}-\mu^2} } \rd r,
    \eeaa
are some elliptic-like integrals. It seems that these cannot be evaluated using the usual elliptic functions.

\bu It is a conventional fact that both $\tpm$ and $\spm$ are analytical functions of $\mu \in(0,1)$. Moreover, whenever $p\ne 2$, it can be seen from Lemma \ref{TSp} that $\tpm$ and $\spm$ are non-constant functions of $\mu \in(0,1)$. For the graphs of these functions, see Figure \ref{fig-TSp}.
\qed
}
\end{ex}

\begin{figure}[httb]
\bc
\resizebox{7cm}{6cm}{\includegraphics{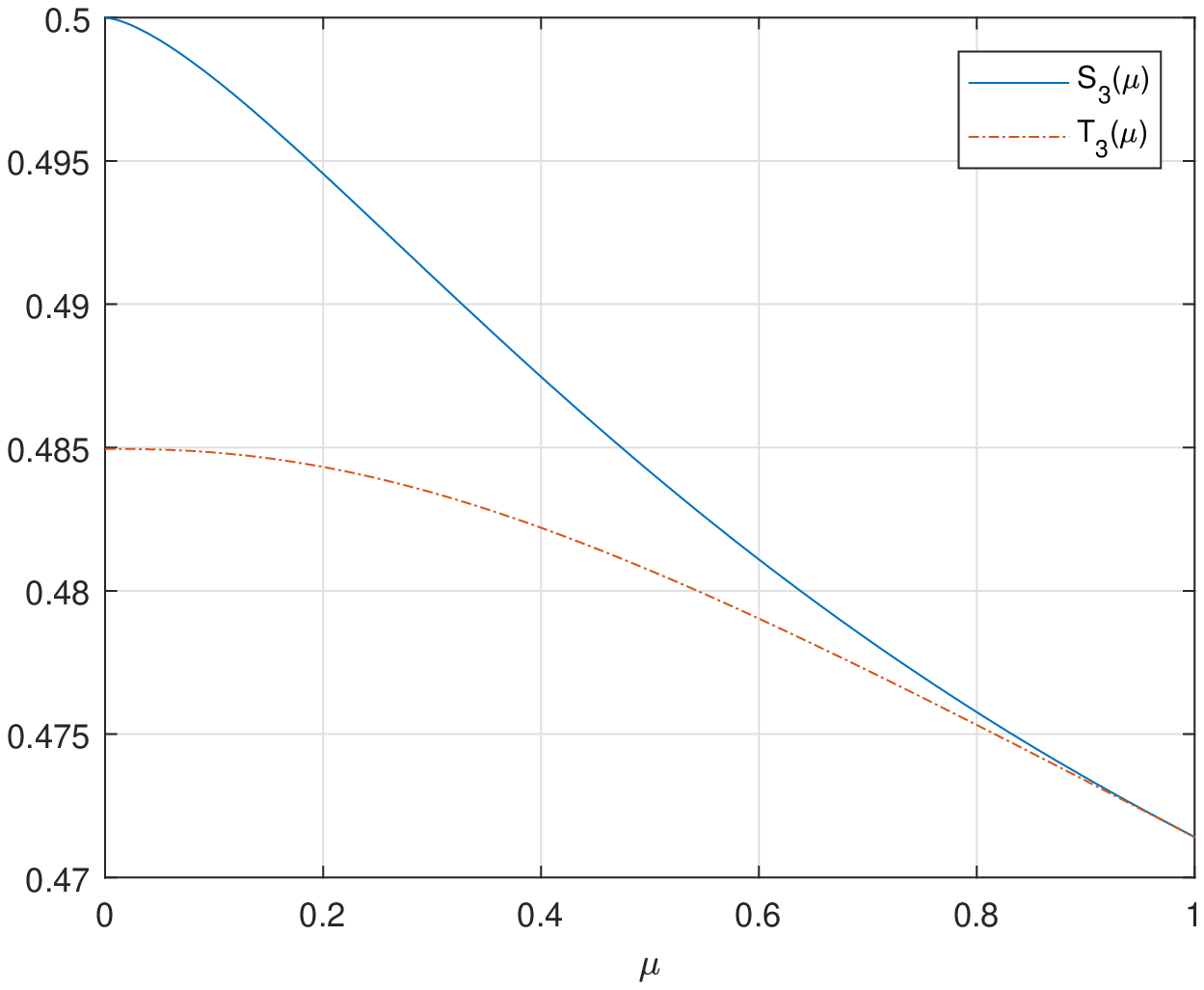}}\hspace{1cm} \resizebox{7cm}{6cm}{\includegraphics{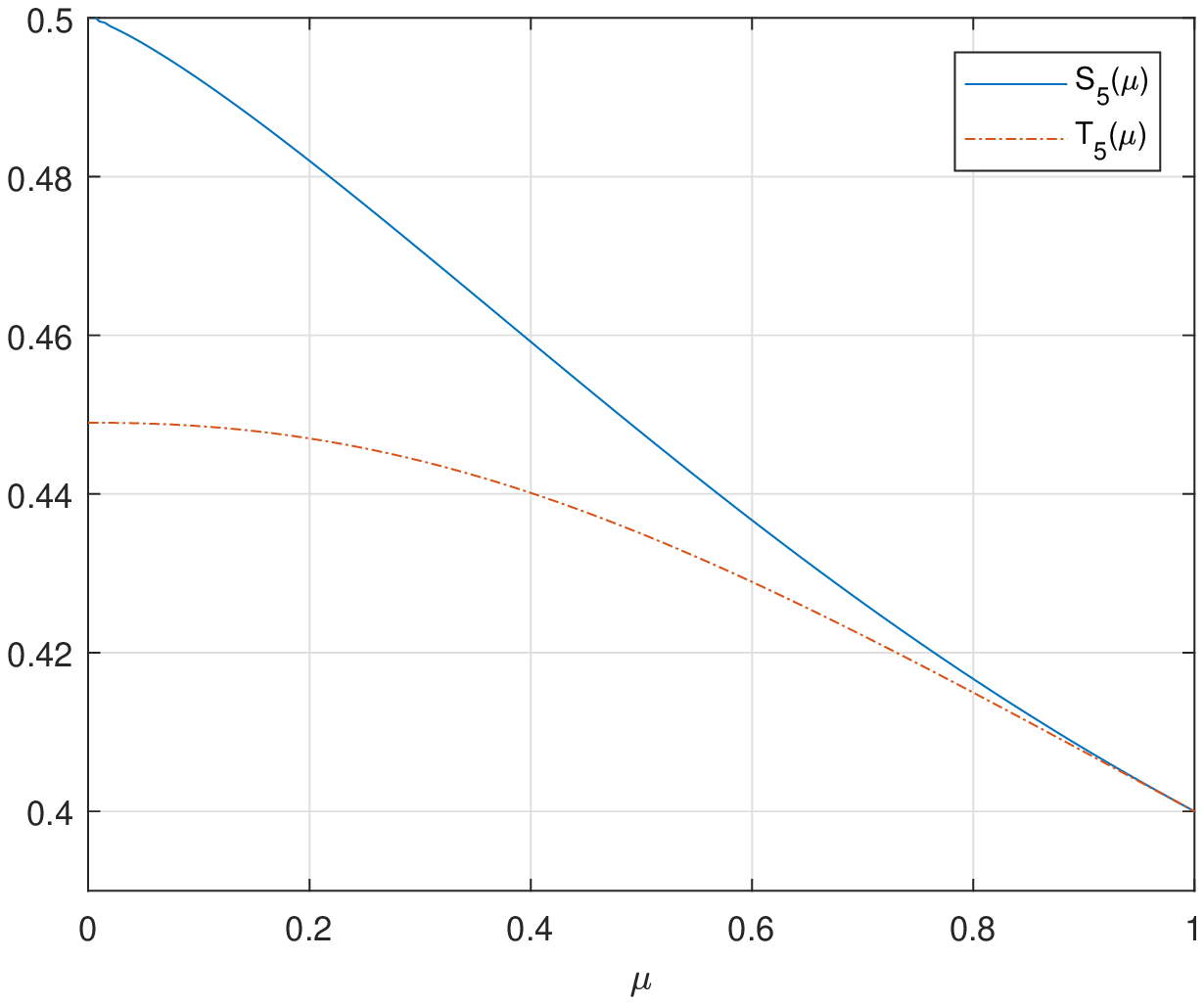}}
\ec
\caption{Functions $\spm$ and $\tpm$, where $p=3$ (left) and $p=5$ (right). }
\lb{fig-TSp}
\end{figure}

Recall that $p$ and $q$ are conjugate exponents. One has the following symmetries of the functions $\tpm$ and $\spm$ in exponents.

\bb{lem} \lb{symm} Whenever $p,\ q\in(1,\oo)$ are conjugate exponents, one has
    \be \lb{TSp1}
    \tpm \equiv T_q(\mu)\andq \spm \equiv S_q(\mu).
    \ee
\end{lem}

\Proof From the defining equalities \x{Fr} and \x{Qr}, one has $F_p(r) \equiv F_q(1-r)$ and $Q_p(r) \equiv Q_q(1-r)$. By \x{Gmu}, $r_{\pm,p}(\mu) \equiv 1- r_{\mp,q}(\mu)$. Hence \x{TSp1} follows immediately from equalities \x{Tmu} and \x{Smu} for $\tpm$ and $\spm$.   \qed


\section{Construction of Periodic Eigenvalues Using the Reduced Dynamics} 
\setcounter{equation}{0} \lb{4th}

This is the main part of the paper. We will first use the dynamics of the reduced system \x{S2d}, i.e., non-constant periodic orbits, to give a complete construction for all possible $1$-periodic eigenvalues of the vectorial $p$-Laplacian for $d=2$. In the second part of this section, the main Theorem \ref{main1} will be proved.

\subsection{Constructions of periodic eigenvalues}  \lb{sect42}

\ifl
Before going the details, let us give another explanation to the eigenvalue sequence $\tlan=2n\pi$, $n\in \N$.

    \bb{ex} \lb{exam2}
{\rm We have known from Proposition \ref{prop1} that, when the scaling momentum $\mu=1$, problem \x{SS}-\x{H11}-\x{P1} of Lemma \ref{spic} leads to eigenvalues $\tlan$$, n\in \N$. We can use the dynamics of the reduced system \x{S2d} to give another explanation. Note that $\mu=1$ corresponds to the equilibrium $C_1$, i.e., to the constant solution
    \(
    (r(t),\th(t))\equiv (1/p,\pi/2).
    \)
In this case, Eq. \x{vpeq} is simply
    \(
    \dot \vp(t) \equiv 1.
    \)
Hence
    \[
    \vp(t) \equiv t+\vp_0, \andq \psi(t) \equiv t+\vp_0+\pi/2.
    \]
By \x{PC}, we have the scaling solutions
    \be \lb{refs}\bb{split}
    \bfx(t) &= e^{i\vp(t)} \equiv e^{i(t+\vp_0)},\andq
    \bfy(t) = (q(1- 1/p))^{1/q} e^{i\psi(t)} \equiv  e^{i(t+\vp_0+\pi/2)}.
    \end{split}
    \ee
Thus \x{P1} implies that $\la=2n\pi$, the eigenvalues $\tlan$ in \x{sp5}. In fact, \x{refs} is just the scaling eigenfunctions for all eigenfunctions $\tl \bfe_p(t;\tlan)$ in \x{sp5}. See \x{sefs0}.\qed
}
    \end{ex}

Now we extend the idea of Example \ref{exam2} to eigenvalues with general scaling momenta.
\fi

Now we establish the relation between the spectral set $\si_p^*$ and the dynamics of the reduced system \x{S2d}. For this purpose, we need only to construct the following subset of the spectral set
    \[
    \si_p^{**}:=\si_p^*\setminus \{\lan,\ \tlan: n \in \N\}.
    \]


    \bb{lem} \lb{rel2}
(1) Let $\la\in \si_p^{**}$. Then there are $\mu\in(0,1)$ and $\ell, \ m\in \N$ such that
    \be \lb{Eq0}
    \la= 2 m \pi T(\mu), \andq m S(\mu)= \ell.
    \ee

(2) Conversely, suppose that $\la$ satisfies \x{Eq0} for some $\mu\in(0,1)$ and $\ell, \ m\in \N$. Then $\la\in \si_p^*$.
    \end{lem}

\Proof We first prove the necessity. Let $\la\in \si_p^{**}$. By Lemma \ref{spic}, $\la$ admits a scaling eigenfunction $\bxy\in \Ga_{1,1,\mu}$ with some scaling momentum $\mu\in(0,1)$ such that $\la$ is a period of $\bxy$. See \x{P1}. Hence $\|\bfx(t)\|$ is also $\la$-periodic. Notice that $\bxy$ corresponds to the $2\pi T$-periodic orbit $C_\mu=\{(r(t),\th(t))\}\subset \D_0^*$ of \x{S2d} of the minimal period $2\pi T$. By \x{PC}, one has
    \[
    \|\bfx(t)\| \equiv (pr(t))^{1/p}.
    \]
By comparing their periods, we conclude that there is some $m\in \N$ such that $\la=2m \pi T$, the first equality of \x{Eq0}.

For the second one, by \x{PC} again,
    \[
    e^{i\vp(t)} =(pr(t))^{-1/p} \bfx(t)
    \]
has the period  $\la=2m \pi T$. Hence there is some $\ell\in \N$ such that $\vp(\la) -\vp(0)= 2\ell \pi$. By using the equality \x{vpt} for the function $S$, we obtain
    \[
    2\ell \pi = \vp(\la) -\vp(0)=\vp(2m \pi T) -\vp(0) =m(\vp(2\pi T)-\vp(0)) = 2 m \pi S,
    \]
the second equality of \x{Eq0}.

Conversely, if \x{Eq0} is satisfied by some $\mu\in(0,1)$ and $\ell, \ m\in \N$, it is necessary that $\la$ is an eigenvalue with some eigenfunction of the scaling momentum $\mu \in(0,1)$. By checking the proof above, one can see that $\la$ is an eigenvalue. However, due to the extremal cases in Example \ref{se-lan}, we can only assert that $\la\in \si_p^*$.
\qed

In the sequel, we are mainly interested in the case $p\ne 2$. However, the following general construction for eigenvalues also includes this trivial case. The two equations in \x{Eq0} of Lemma \ref{rel2} can be used to construct all periodic eigenvalues $\la\in \si_p^{**}$ of problem \x{R}-\x{bc}. For convenience, let us consider $\spm$ as a function of the scaling momenta $\mu\in (0,1)$, whose range is denoted by
    \[
    \mcs_p:= \z\{\spm: \mu\in(0,1) \y\}\subset (0,+\oo).
    \]

For any rational number $\lm\in \mcs_p^*$, where $\ell, \ m\in \N$ and the irreducibility of $\lm$ is not assumed at this moment, there is at least one scaling momentum $\mu_\plm\in(0,1)$ such that
    \be \lb{Eq1}
    S_p(\mu_\plm) = \lm\in \mcs_p^*.
    \ee
Then we can define a number
    \be \lb{Eq2}
    \pilm:= m\pi T_p(\mu_\plm)\equiv \ell \pi\f{ T_p(\mu_\plm)}{S_p(\mu_\plm)}>0.
    \ee
Lemma \ref{rel2} asserts that $2\pilm\in \si_p^*$. Moreover, such an eigenvalue $2\pilm$ will admit some eigenfunctions of the scaling momentum being $\mu_\plm$.

    \bb{rem}\lb{ev-const}
The idea of the construction of $2\pilm\in \si_p^*$ is very simple. Whenever $p\ne 2$, by assuming that $\spm$ is strictly decreasing in $\mu\in(0,1)$, which will be conjectured in Section \S \ref{sec6},  one has
    \be \lb{evs-lm}
    \mu_\plm \equiv S_p^{-1}(\lm), \andq \pilm \equiv m\pi T_p\z( S_p^{-1}(\lm)\y).
    \ee
Then $2\pilm\in \si_p^*$ is an eigenvalue, while the $\mu_\plm$ of \x{evs-lm} is the scaling momentum of some eigenfunctions associated with $2\pilm$. However, though we are not able to prove that $\spm$ is strictly decreasing in the whole interval $(0,1)$, this is true when $\mu$ is near $1$. Actually, such a local monotonicity is very important in the proof of Theorem \ref{main}.
    \end{rem}

    \bb{ex} \lb{2lap}
{\rm Let us consider the Laplacian, i.e. $p=2$. One has $S_2=\{1/2\}$. Hence Eq. \x{Eq1} is only fulfilled by choosing $\ell=n\in \N$ and $m=2n$. By \x{Eq2} and \x{T2mu}, \x{S2mu},
    \[
    \pi_{2,n/2n} = 2n \pi \d \f{1}{2} \equiv n \pi\qqf n\in \N.
    \]
Hence the resulted eigenvalues are $2n \pi$. These are $1$-periodic eigenvalues for the Laplacian, but are not new to us.\qed}
    \end{ex}

On the other hand, it is trivial from \x{Eq1} and \x{Eq2} that whenever $\lm\in \mcs_p$, by setting $\ell$ and $m$ as $n \ell$ and $n m$, one has
    \be \lb{Re1}
    \mu_{p,n \ell/n m} \equiv \mu_\plm, \andq \pi_{p,n \ell/n m} \equiv n\pi_\plm\qqf n\in \N.
    \ee
By Lemma \ref{rel2} again, $2n\pilm\in \si_p^*$ are eigenvalues, while the scaling momenta still are $\mu_\plm$. In order to eliminate the dilation relations \x{Re1} for eigenvalues and scaling momenta, in the sequel, let us only choose  in \x{Eq1}-\x{Eq2} the {\it irreducible} rational numbers $\lm$ from the following subset
    \be \lb{sps}
    \mcs_p^*:= \z\{ \lm\in \mcs_p\cap \Q: \mbox{$\ell, \ m\in \N$ are co-prime}\y\}.
    \ee
For $p=2$, we have seen from \x{S2mu} that $\mcs_2^*=\{1/2\}$, a single point. Moreover, whenever $p\ne 2$, $\mcs_p$ is a nontrivial interval and $\mcs_p^*$ always contains infinitely many different rational numbers. See Lemma \ref{TSp} and Remark \ref{int} below.

For any $\lm\in \mcs_p^*$, one has constructed $\pilm$ so that
    \be \lb{Eq3}
    \La_\plm:=\z\{ \tlann:=2 n\pilm: n\in \N\y\}\subset \si_p^*.
    \ee
Due to the extremal case as in Example \ref{2lap} and some additional difficulties for the $p$-Laplacian, the following constructions for the spectral sets $\si_p^{**}$ and $\si_p^*$ are stated in a little bit complicated way.
\ifl For later convenience, we also denote
    \[
    \bb{split}
    \pi_{p,0}:= \pi_p, \qq     & \La_{p,0}:=\z\{ 2 n\pi_p: n\in \N\y\}\\
    \pi_{p,\oo}:= \pi, \qq     & \La_{p,\oo}:=\z\{ 2 n\pi: n\in \N\y\}.
    \end{split}
    \]
\fi

    \bb{thm} \lb{main1} Let $p\in (1,\oo)$. Then
    \bea
    \si_p^{**} \AND{\subset} \z\{\La_\plm: \lm\in \mcs_p^*\y\}\nn
    \\
    \EE \z\{2 n\pi_*: n\in \N,\ \mbox{ and $\pi_*$ stands for either $\pi_p$ or $\pi$ or $\pilm$ with $\lm \in \mcs_p^*$} \y\},\nn \\
    \lb{sip2}
    \si_p^{*} \EE \z\{\lan, \ \tlan:n\in \N \y\}\cup \z\{\La_\plm: \lm\in \mcs_p^*\y\}\nn\\
    \EE \z\{2 n\pi_*: n\in \N, \ \mbox{ and $\pi_*$ stands for either $\pi_p$ or $\pi$ or $\pilm$ with $\lm \in \mcs_p^*$} \y\}.
    \eea
    \end{thm}

\ifl
Here one can notice for the case $p=2$, the inclusion in \x{sip1} cannot be improved as equality. See Example \ref{2lap}.
while \x{sip2} can be written as an equality by understanding
    \[
    \La_{2,0}\equiv \La_{2,\oo}\equiv \La_{2,\lm}\equiv \z\{ 2 n\pi: n\in \N\y\}\qqf \lm \in \si_2^*,
    \]
because $\si_2^*=\{1/2\}$ contains only one irreducible number $1/2$. See \x{pi0oo}.

{\bf Proof of Theorem \ref{main1}} Given any irreducible rational number $\lm$, let $\mu=\mu_\plm\in(0,1)$ be as in \x{per3}. Then for any $n\in \N$, one also has
    \[
    S_p(\mu_\plm) =2 \pi \f{\ell n}{m n},
    \]
which corresponds to \x{per1}, where $\ell$ and $m$ are replaced by $\ell n$ and $m n$ respectively. By using the corresponding equality \x{per1}, we know that
    \[
    \tlann= m n T_p(\mu_\plm)=2 n\pilm
    \]
is a $1$-periodic eigenvalue of problem \x{R} whose eigenfunctions correspond to non-constant periodic orbits of system \x{S2d}.

By checking the constructions above carefully and combining with Theorems \ref{prop1}, we know that the converse of the above is also true.
\qed
\fi

It follows from the symmetries \x{TSp1} in Lemma \ref{symm} and Theorem \ref{main1} the following results.

    \bb{lem} \lb{ev-pq}
Whenever $p,\ q\in(1,\oo)$ are conjugate exponents, one has
        \[
        \si_p^* =\si_q^*.
        \]
Going to the original spectral sets $\Si_p$ of problem \x{pd}-\x{bc}, the symmetry in the conjugate exponents can be informally written as
    \[
    \z(\Si_p\y)^{1/p} \equiv \z(\Si_q\y)^{1/q} \qq \mbox{for conjugate exponents $p, \ q$}.
    \]
    \end{lem}

Results \x{Eq3} and \x{sip2} assert that each sequence takes the form $\{2n\pi_{p,*}: n\in \N\}$. However, once these are used to
obtain new eigenvalue sequences other than the known ones, there arises a very delicate problem, as explained in the remarks after Theorem \ref{main}. To address this issue, let us introduce the following notion.

    \bb{defn} \lb{nind}
{\rm We say that a set $W \subset (0,\oo)$ is {\it integer-independent}, if
    \be \lb{ww}
    \f{w}{w'} \not \in \N \q \mbox{for any different } w, \ w'\in W.
    \ee
}
    \end{defn}

For example, when $p\ne 2$, we have known that $\{\pi_p,\ \pi\}$ is integer-independent. We remark that, when $W$ consists of $k<\oo$ numbers, the integer-independence \x{ww} includes $N_k:=k^2-k$ inequalities.

\ifl
we must consider the case $p\ne 2$. Suppose that some irreducible number $\lm\in \mcs_p^*$ is used to construct an eigenvalue sequence $\La_\plm$ as in \x{Eq3}. If $\pi_\plm$ satisfies
    \be \lb{pis}
    \mbox{either} \q \f{\pi_\plm}{\pi_{p,0}} =k \in \N, \q \mbox{or}\q \f{\pi_\plm}{\pi_{p,\oo}} =k \in \N,
    \ee
say the former case, then the eigenvalue sequence
    \[
    \La_\plm= \z\{2n \pi_\plm: n\in \N \y\}= \z\{2n k\pi_{p,0}: n\in \N \y\}
    \]
is the same as $\La_{p,0}$ itself if $k=1$, and is a proper subsequence of $\La_{p,0}$ if $k>1$. Hence $\La_\plm$ has actually not provided any new eigenvalue different from $\La_{p,0}\cup \La_{p,\oo}$. In order that \x{Eq3} can really yield a new eigenvalue sequence of problem \x{R}-\x{bc}, we need to exclude the cases in \x{pis}. 

    \bb{prop} \lb{pev-n}
Let $p\ne 2$ and $\pilm$ be as in \x{Eq2}. Then $\La_\plm$ of \x{Eq3} is a sequence of $1$-periodic eigenvalues of \x{R}-\x{bc} so that it is different from the known sequences $\La_{p,0}$ and $\La_{p,\oo}$ iff the following $4$ requirements are fulfilled
    \be \lb{new1}
    \f{\pilm}{\pi}\not \in \N,\qq \f{\pilm}{\pi_p}\not \in \N, \qq \f{\pi}{\pilm}\not \in \N,\qq \f{\pi_p}{\pilm}\not \in \N.
    \ee
\end{prop}

The latter two requirements of \x{new1} have excluded the case that the original sequences $\{\lan\}$ and $\{\tlan\}$ become a subsequence of the newly constructed sequence $\La_\plm$.

A more delicate problem in the obtention of new eigenvalues is as follows. Suppose that we have added to $\si_p^*$ a new sequence $\La_\plm$ as in Proposition \ref{pev-n}. In order that $\La_{p,\ell'/m'}$, $\ell'/m'\in \mcs_p^*$ to be a new sequence of eigenvalues different from $\{\lan\}$ and $\{\tlan\}$ and different from $\La_\plm$ as well, besides comparing $\pi_{p,\ell'/m'}$ with $\pi_p, \ \pi$ like \x{new1}, one needs also to verify
    \be \lb{new2}
    \f{\pi_{p,\ell'/m'}}{\pilm}\not \in \N, \andq \f{\pilm}{\pi_{p,\ell'/m'}}\not \in \N.
    \ee
The more sequences we have at hand, the more requirements like \x{new1} and \x{new2} should be verified.
\fi

\subsection{Proof of the main theorem}\lb{main-p}

Due to Theorem \ref{main1} and Lemma \ref{ev-pq}, we can restrict the exponent $p$ to the interval $(1,2)$ or to the interval $(2,\oo)$. In the sequel, let $p>2$ be fixed. For evidence, in the sequel, we will not use the conjugate exponent $q$ no longer and write down all results using the exponent $p$ itself.

Firstly we give some asymptotic properties of the functions $T(\mu)=\tpm$ and $S(\mu)= \spm$ near $\mu=1$.  Suggested by \x{Eq2}, we also need to study the ratio function
    \be \lb{Umu}
    U(\mu)= \upm:= \f{\tpm}{\spm}, \qq \mu \in(0,1).
    \ee
The following proofs are borrowed from the techniques for bifurcation theory and limit cycles. For simplicity, we use $C_i=C_i(p)$ to denote constants depending only on $p$, most of them are positive.

\begin{lem}\label{TSp}
    When $\mu\rightarrow 1^-$, there have
    \begin{eqnarray}\label{TS1}
    \displaystyle T(\mu) \EQ C_1-C_2(\mu-1)+ C_3(\mu-1)^2+o((\mu-1)^2), \\
    \label{TS2}
    \displaystyle S(\mu) \EQ C_1-C_2(\mu-1)+ C_4(\mu-1)^2+o((\mu-1)^2),
    \end{eqnarray}
where
    \[
    \bb{split}
    C_1 &:= \frac{\sqrt{p-1}}{p}>0,\qq \qq C_2:=\frac{(p-2)^2}{12 p\sqrt{p-1}}>0, \\
    C_3 &:= \frac{(p-2)^4}{576 p \z(\sqrt{p-1}\y)^3}>0,\qq C_4:=\frac{(p-2)^2(p^2-20p-20)}{576 p\z(\sqrt{p-1}\y)^3}.
    \end{split}
    \]
\end{lem}

\Proof Let $y=\sqrt{Q^2(r)-\mu^2}$ and consider closed curves
    $$
    E_\mu:\q \{(r,y): y^2+1-Q^2(r)=1-\mu^2\}, \qq \mu\in(0,1).
    $$
When $\mu\rightarrow 1^-$, $E_\mu$ tends to the point $(1/p, 0)$. Notice that
    $$1-Q^2(r)=\frac{p^2}{p-1}\z(r-\frac{1}{p}\y)^2+o\z(\z(r-\frac{1}{p}\y)^2\y).$$
We can find a transformation
    \be \lb{ru}
    r-\frac{1}{p}=\frac { \sqrt{p-1}}{p}u+\frac{p-2}{3p}{u}^{2}+\frac{(p-2)^2}{36p \sqrt{p-1}}{u}^{3}
+C_5{u}^{4}+C_6{u}^{5}+o(u^5),
    \ee
such that
    $$1-Q^2(r)=u^2,$$
where
    $$
    C_5=-\frac{(p-2)(p^2-22p+22)}{ 270(p-1) p}, \qquad C_6={\frac {(p-2)^2(p^2+68p-68)}{4320\,p
\z(\sqrt{p-1}\y)^3}}.
    $$
That is, the curve $E_\mu$ becomes a circle
    $$
    E_\mu: \q \{(u,y):y^2+u^2=1-\mu^2\}.
    $$

Using the transformation \x{ru} above, the integrals $T(\mu)$ and $S(\mu)$ of \x{Tmu} and \x{Smu} are now
    \[
T(\mu)=\f{1}{\pi} \int_{-\sqrt{1-\mu^2}}^{\sqrt{1-\mu^2}}\frac{f_1(u)}{y}\du, \andq
S(\mu)= \f{\mu}{\pi} \int_{-\sqrt{1-\mu^2}}^{\sqrt{1-\mu^2}}\frac{f_2(u)}{y}\du,
    \]
where the functions $f_1$ and $f_2$ can be calculated explicitly. When $\mu\rightarrow 1^-$, by finding the Taylor series of $f_1$ and $f_2$, we can obtain results \x{TS1} and \x{TS2}.   \qed

    \bb{rem} \lb{int}
For any $p\ne 2$, one has $C_2\ne 0$. Thus we see from \x{TS1} and \x{TS2} that both $\tpm$ and $\spm$ are non-constant functions.  In particular, the range $\mcs_p$ of $S_p(\d)$ is a non-trivial interval, and the set $\mcs_p^*$ we have used in \x{sps} is necessarily an infinitely many set.
    \end{rem}

\medskip
Due to \x{TS1} and \x{TS2}, it is conventional to define
    \be \lb{TS11}
    T(1)=S(1):= C_1=\frac{\sqrt{p-1}}{p}.
    \ee
From Lemma \ref{TSp}, some observations on $T(\mu)$ and $S(\mu)$ are stated as the following lemma.

    \bb{lem} \lb{tan}
(1) By extending $T(\mu)$ and $S(\mu)$ to $\mu=1$ as in \x{TS11}, the functions $T(\mu)$ and $S(\mu)$ have the same tangent line at $\mu=1$. Moreover, one has the following asymptotic expansion for $U(\mu)$
    \begin{equation}\lb{TSW}
    U'(\mu)=-C_7(\mu-1)+o((\mu-1)), \qquad \mu\rightarrow 1^-,
    \end{equation}
where
    \[
    C_7:=\frac{(p-2)^2}{12(p-1)}>0.
    \]

(2) Consequently, there is some $\delta_0\in (0, 1)$ such that both $T(\mu)$ and $S(\mu)$ are strictly decreasing on $(1-\delta_0, 1)$, and
    \be \lb{est12}
    |S'(\mu)|\geq C_2/2,\andq |U'(\mu)|\leq 2 C_7(1-\mu) \qqf \mu\in(1-\delta_0, 1).
    \ee
    \end{lem}

In fact, \x{TS1} and \x{TS2} imply that $T'(1)=S'(1)= -C_2$. Moreover, from the definition \x{Umu} of $U(\mu)$, the expansion  \x{TSW} can also be deduced. Then estimates \x{est12} are following from \x{TS1}, \x{TS2} and \x{TSW}. For the graph of $\upm$, see Figure \ref{fig-Up}.

\begin{figure}[httb]
\bc
\resizebox{8cm}{7cm}{\includegraphics{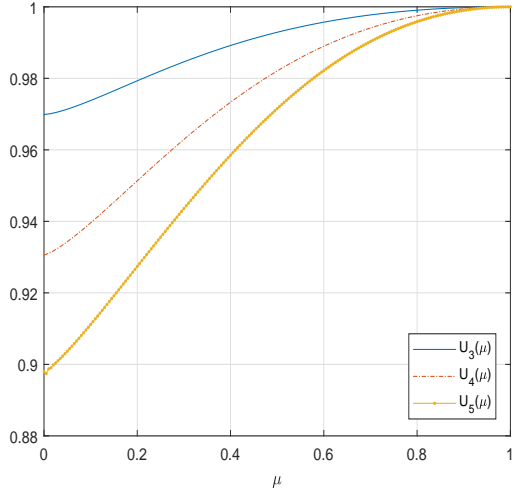}} 
\ec
\caption{Functions $\upm$, where $p=3$ (top), $p=4$ (middle), and $p=5$ (bottom). }
\lb{fig-Up}
\end{figure}

%
\medskip

Now we give the proof of Theorem \ref{main}. 
Note that $S(1)=C_1$ is as in \x{TS11}. For convenience, denote
    \[
    \delta :=S(1-\delta_0)-S(1)>0.
    \]
Since there are only finite positive integers $n$ such that $\frac{1}{n}>C_1$, we take smaller $\da, \ \da_0$ so that
    \be \lb{N0}
    \frac{1}{n}\not\in (C_1, C_1+\delta)\qqf n\in \N.
    \ee

In the sequel, let $n\in \N$ be arbitrarily given. It is possible to find large prime numbers $\ell_n\uto+\oo$ such that
    $$
    \frac{\ell_n}{C_1}- \frac{\ell_n}{C_1+\frac{\delta}{2^n}}\ge n+2.
    $$
Hence the open interval
    $$\z(\frac{\ell_n}{C_1+\frac{\delta}{2^n}}, \frac{\ell_n}{C_1}\y)$$
contains at least $n+1$ integers $m_n,\ m_n+1,\ \dd, \ m_n+n$. In other words, for any $i\in \{0, 1, \dd, n\}$, one has
    \be \lb{lni0}
    \frac{\ell_n}{m_n+i}\in \z(C_1, C_1+\frac{\delta}{2^n}\y)\subset (C_1, C_1+\delta).
    \ee
Obviously, each $\ell_n/(m_n+i)$ is irreducible. Otherwise, as $\ell_n$ is a prime, we have $m_n+i = k\ell_n$ for some $k\in\N$. Consequently,
    \[
    \f{1}{k} \in\z(C_1, C_1+\frac{\delta}{2^n}\y)\subset (C_1, C_1+\delta),
    \]
a contradiction to the construction \x{N0} for $\da$. As a consequence of the irreducibility, one sees that
    \be \lb{lni}
    \ell_n/(m_n+i), \q i=0,1,\dd,n
    \ee
must be different when $n$ is given.

Due to the construction of $C_1$ and $\da$, it follows from \x{lni0} and \x{lni} that
we can introduce
    \be \lb{muni}
    \mu_{n, i}\in (1-\delta_0, 1)\q \mbox{such that}\q S( \mu_{n, i})=\ell_n/(m_n+i), \q i=0,1,\dd,n.
    \ee
Define then
    \be \lb{pini}
    \vpi_{n, i}:= 
    \ell_n \pi U(\mu_{n,i})\in\si_p^*, \q i=0,1,\dd,n.
    \ee
From the properties of $S(\mu)$ and $U(\mu)$, one sees that whenever $n$ is fixed, $\vpi_{n, i}$, $i=0,1,\dd,n$ are different. The bi-sequence \x{pini} is the candidate for our construction of new eigenvalues.

We are going to prove that $\si_p^*$ contains infinitely many different sequences of eigenvalues stated in Theorem \ref{main}. Otherwise, let us assume that 
    \[
    \si_p^*= \z\{2n\vpi_1:n\in \N\y\}\cup \z\{2n\vpi_2: n\in \N\y\}\cup \dd\cup \z\{2n\vpi_k:n\in \N\y\},
    \]
where $k\in \N$ is finite and $0<\vpi_1<\dd<\vpi_k$. Hence the integer-independence of $W_k:=\{\vpi_1, \dd, \vpi_k\}$ can be assumed.

Let $n\geq k$ be arbitrarily given. At first we consider the first $(k+1)$ numbers $\vpi_{n, 0}, \ \vpi_{n, 1},\ \ldots,\ \vpi_{n, k}$ in \x{pini}. Because of \x{muni} and \x{pini}, in order to complete the proof of Theorem \ref{main}, it suffices to prove the following claim: Whenever $n$ is large enough, there must be some $i_0\in\{0,1,\dd,k\}$, depending on $n$ as well, such that
    \be \lb{claim}
    W_{k+1,n} := \z\{ \vpi_1, \dd, \vpi_k, \, \vpi_{n,i_0}\y\}\q \mbox{is integer-independent}.
    \ee
By the explanation following Definition \ref{nind}, claim \x{claim} is equivalent to $N_{k+1}-N_k =2 k$ inequalities, which can be classified into the following two classes
    \bea \lb{cl1}
    \f{\vpi_j}{\vpi_{n,i_0}}\AND{ \not \in } \N \qqf j=1,\dd,k,\\
    \lb{cl2}
    \f{\vpi_{n,i_0}}{\vpi_j}\AND{ \not \in } \N \qqf j=1,\dd,k.
    \eea

The requirement in \x{cl1} is simple. Let us emphasize the constructions in \x{muni} and \x{pini} as follows. For any $i=0,1,\dd,k$, one has
    \be \lb{muni1}
    \mu_{n,i} \in (1-\da_0,1), \andq U(\mu_{n,i}) \mbox{ is close to $1$}.
    \ee
Therefore we are able to obtain
    \be \lb{cl11}
    \vpi_{n,i} = \ell_n \pi U(\mu_{n,i}) \sim \ell_n \pi.
    \ee
As the primes $\ell_n$ are large when $n\to \oo$, we know that, if $n$ is large enough,
    \[
    \f{\vpi_j}{\vpi_{n,i}}\sim 0  \qqf j=1,\dd,k,\ i=0,1,\dd,k.
    \]
Hence \x{cl1} is verified for all $i_0=0,1,\dd,k$ as long as $n$ is large.

To verify the requirement in \x{cl2}, let us assume that \x{cl2} is false. That is, by fixing $n\gg k$, we have that for any $i=0,1,\dd,k$, $\vpi_{n,i}$ is an integer multiple of some $\vpi_j$, $1\leq j\leq k$. Hence we can find some $j_0$ and $i_1$, $i_2$ such that $1\le j_0\le k$, $0\leq i_1<i_2\leq k$, and
    $$
    \frac{\vpi_{n, i_1}}{\vpi_{j_0}}=n_1\in \mathbb N, \qquad \frac{\vpi_{n, i_2}}{\vpi_{j_0}}=n_2\in \N,
    $$
i.e.
    $$\ell_n\pi U(\mu_{n, i_1})=n_1\vpi_{j_0}, \qquad \ell_n\pi U(\mu_{n, i_2})=n_2\vpi_{j_0}.$$
Obviously $n_1\ne n_2$. Thus
    \begin{equation}\label{ine1}
    0<\vpi_1\leq \left |n_2-n_1 \right |\vpi_{j_0}=\ell_n\pi \left |U(\mu_{n, i_2})-U(\mu_{n, i_1})\right |.
    \end{equation}

Now we give some necessary estimates. By \x{lni0}, one has
    $$\frac{\ell_n}{m_n+i_1}, \ \frac{\ell_n}{m_n+i_2}\in \z(C_1, C_1+\frac{\delta}{2^n}\y).$$
We have then from \x{est12}
    \[ \bb{split}
    \f{C_2}{2} |\mu_{n, i_1}-1| \leq &\left |S(\mu_{n, i_1})-S(1)\right |=\frac{\ell_n}{m_n+i_1}-C_1 \leq \frac{\delta}{2^n}, \\
    \f{C_2}{2}|\mu_{n, i_1}-\mu_{n, i_2}|\leq &\left |S(\mu_{n, i_1})-S(\mu_{n, i_2})\right |=\frac{\ell_n}{m_n+i_1}-\frac{\ell_n}{m_n+i_2}\leq \frac{\ell_nk}{m_n^2}.
 \end{split}
    \]
Hence
    \begin{equation}\lb{est3}
|\mu_{n, i_1}-1|\leq C_8\frac{\delta}{2^n}, \qquad  |\mu_{n, i_1}-\mu_{n, i_2}|\leq C_8\frac{\ell_nk}{m_n^2},
    \end{equation}
where $C_8:=2/C_2>0$. With these estimates \x{est3} at hand, we have again from \x{est12} that
    \[
\left |U(\mu_{n, i_2})-U(\mu_{n, i_1})\right |\leq 2 C_7(1-\mu_{n, i_1})|\mu_{n, i_1}-\mu_{n, i_2}|\leq
C_9 \delta k\frac{ \ell_n }{2^n m_n^2},
    \]
where $C_9:=2 C_7 C_8^2>0$. Combining with \eqref{ine1}, we obtain the following inequality
    \[
0<\vpi_1 \leq \frac{\pi C_9\delta k}{2^n}\z(\frac{\ell_n}{m_n}\y)^2\leq \frac{\pi C_9\delta k}{2^n}(C_1+\delta)^2,
    \]
which is impossible if $n\gg k$ is large enough. Such a contradiction has verified \x{cl2} for all $n$ large enough.

As mentioned before, the proof to Theorem \ref{main} has been completed. \qed

The last step of the proofs presented here for Theorem \ref{main} is not constructive. However, the new eigenvalue sequences are constructed using some constant $\da_0$ to control the scaling momenta for the new eigenvalues. Thus \x{muni1} and \x{cl11} can yield some useful information on eigenvalues and scaling momenta constructed in this paper. In fact, as $\da_0>0$ can be taken to be arbitrarily small, all of the new eigenvalue sequences in \x{vpim}--\x{evs-s} of Theorem \ref{main} can admit eigenfunctions of scaling momenta close to the maximal momentum $1$. More precisely, let us state the following results which can suggest further study to the structure of spectral sets $\si_p^*$.

    \bb{cor} \lb{pana1}
Let $p\ne 2$ be given. Then the eigenvalue sequences $\La_m=\{(2n \vpi_m)^p: n\in \N\}$, $m=1,2,\dd$ in \x{vpim}--\x{evs-s} of Theorem \ref{main} can be chosen so that
    \be \lb{nevm}
    \lim_{m\to\oo} \vpi_m =+\oo.
    \ee
Moreover, any eigenvalue from $\La_m$ will admit some eigenfunction of the scaling momentum $\mu_m$ satisfying
    \be \lb{nmm}
    \lim_{m\to\oo} \mu_m =1.
    \ee
    \end{cor}

Here \x{nevm} and \x{nmm} can be deduced from estimates \x{muni1} and \x{cl11} in the proof. Going to the spectral set, these mean that there are many different eigenvalues tending to infinity, and the corresponding scaling momenta are all close to the maximal scaling momentum. In such a sense, the panorama of the structures of the spectral sets and their scaling momenta is far from being understood completely.


\section{Numerical Simulations}
\setcounter{equation}{0} \lb{sec5}

Since the kernel idea of this paper is to use the scaling momenta of the eigenfunctions to understand the structure of eigenvalues, besides some numerical simulations to eigenvalues, we are also giving simulations to eigenfunctions, i.e. the periodic motions themselves. As the Hamiltonian system \x{HS} has degree $2$ of freedom, it is luck that we can draw the motions in a proper way.

We state an asymptotical result without proof. Let $p\in(1,\oo)$ be fixed. We have known that
    \be\lb{TL}
    \lim_{\mu\uparrow 1} \tpm = \lim_{\mu\uparrow 1} \spm= \sqrt{p-1}/p.
    \ee
It can also be proved that
    \be  \lb{SL}
    \lim_{\mu\downarrow 0} \tpm = \pi_p/2\pi,\andq \lim_{\mu\downarrow 0} \spm = 1/2.
    \ee
See Figure \ref{fig-TSp}.

Now we give some numerical simulations to eigenvalues and scaling momenta.

    \bb{ex} \lb{p=3}
{\rm For $p=3$, one has from \x{TL} and \x{SL} that
    \[
    \mcs_3\supset (\sqrt{2}/3,1/2)=\z(0.4714, 0.5\y).
    \]
Then $\mcs_3^*$ contains
    \[
    \lm = 9/19=0.4737.
    \]
One can obtain the numerical results
    \be \lb{mulap3}
    \mu_{3,9/19}=0.8906\andq \pi_{3,9/19}=28.2668
    \ee
Notice that the two ratios are
    \[
    \f{\pi_{3,9/19}}{\pi}=8.9976 <9 < \f{\pi_{3,9/19}}{\pi_3}=9.2769,
    \]
and their reciprocals are $<1$. Hence $\pi_{3,9/19}$ can yield a new eigenvalue sequence.

Going back to the original problem \x{R}-\x{bc}, the least one of this sequence of eigenvalues deduced from \x{mulap3} is
    \[
    \la_{3,9/19}=(2 \pi_{3,9/19})^3=18068.4095,
    \]
with an eigenfunction
    \[
    \bfe_{3,9/19}(t):=\bfe_3(t; \la_{3,9/19})=(x_1(t),x_2(t))
    \]
of the scaling energy $1$ and the scaling momentum $\mu_{3,9/19}$ as in \x{mulap3}.
    \qed
}
    \end{ex}

The periodic motions of $\bfe_{3,9/19}(t)$ in different spaces are plotted in Figure \ref{persolp3}. These are apparently like the Lissoajous figures.

\begin{figure}[httb]
\bc
\resizebox{7cm}{6cm}{\includegraphics{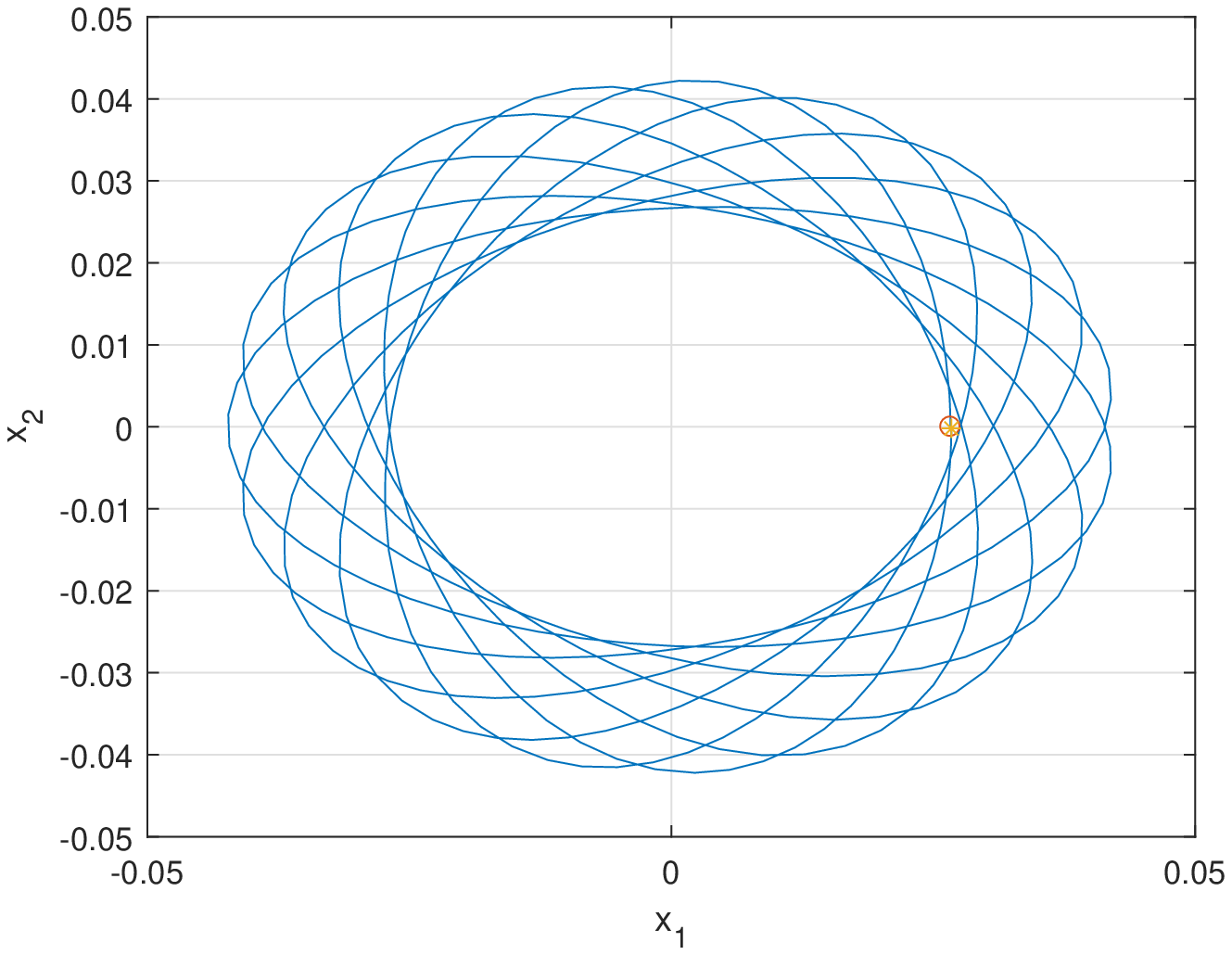}}\hspace{1cm} \resizebox{7cm}{7cm}{\includegraphics{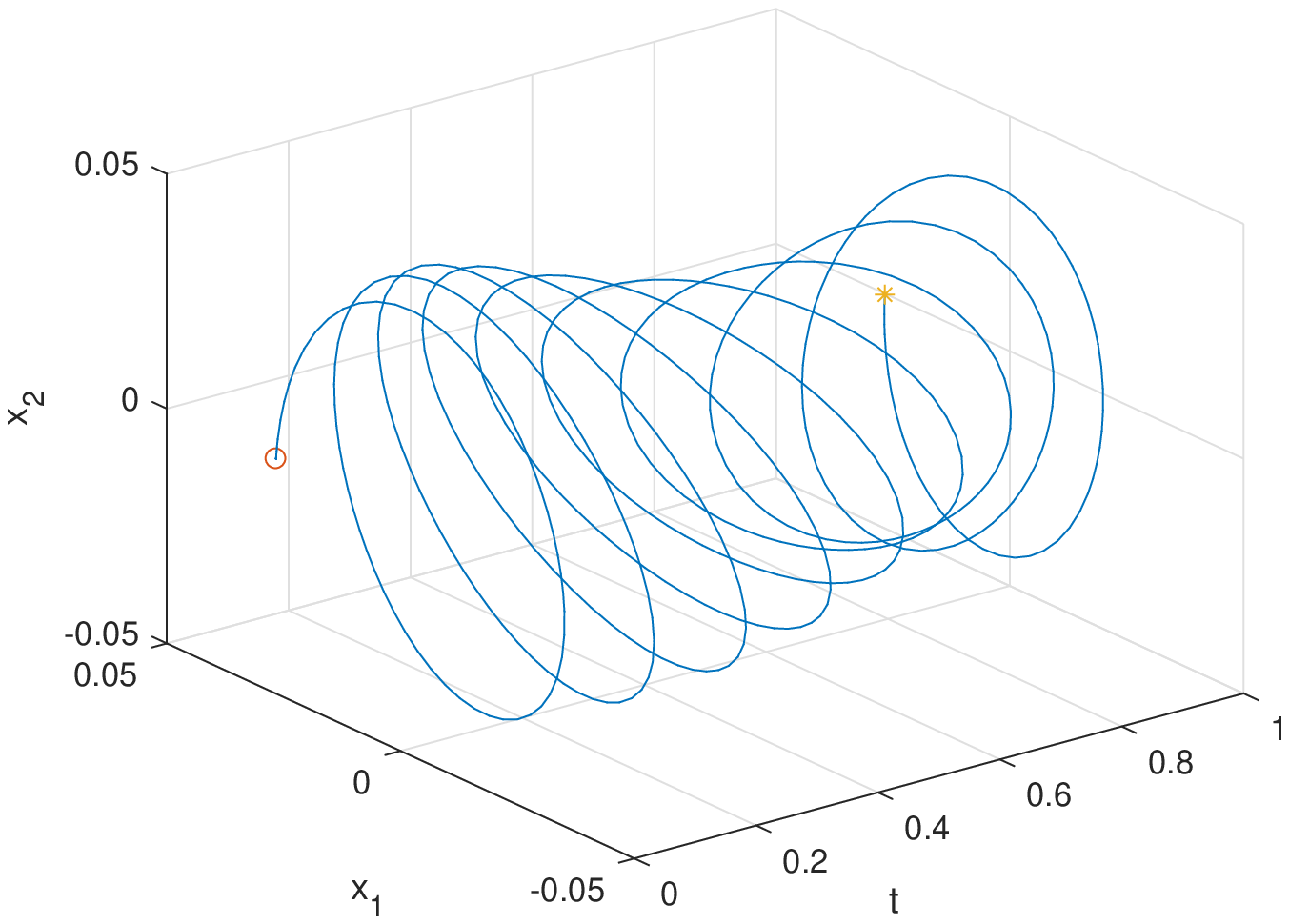}}
\ec
\caption{The periodic orbits (left) and trajectories (right) of the eigenfunction $\bfe_\plm(t)$. Here $p=3$ and $\lm=9/19$.
}
\lb{persolp3}
\end{figure}

\ifl
\begin{figure}[httb]
\bc
\resizebox{7cm}{7cm}{\includegraphics{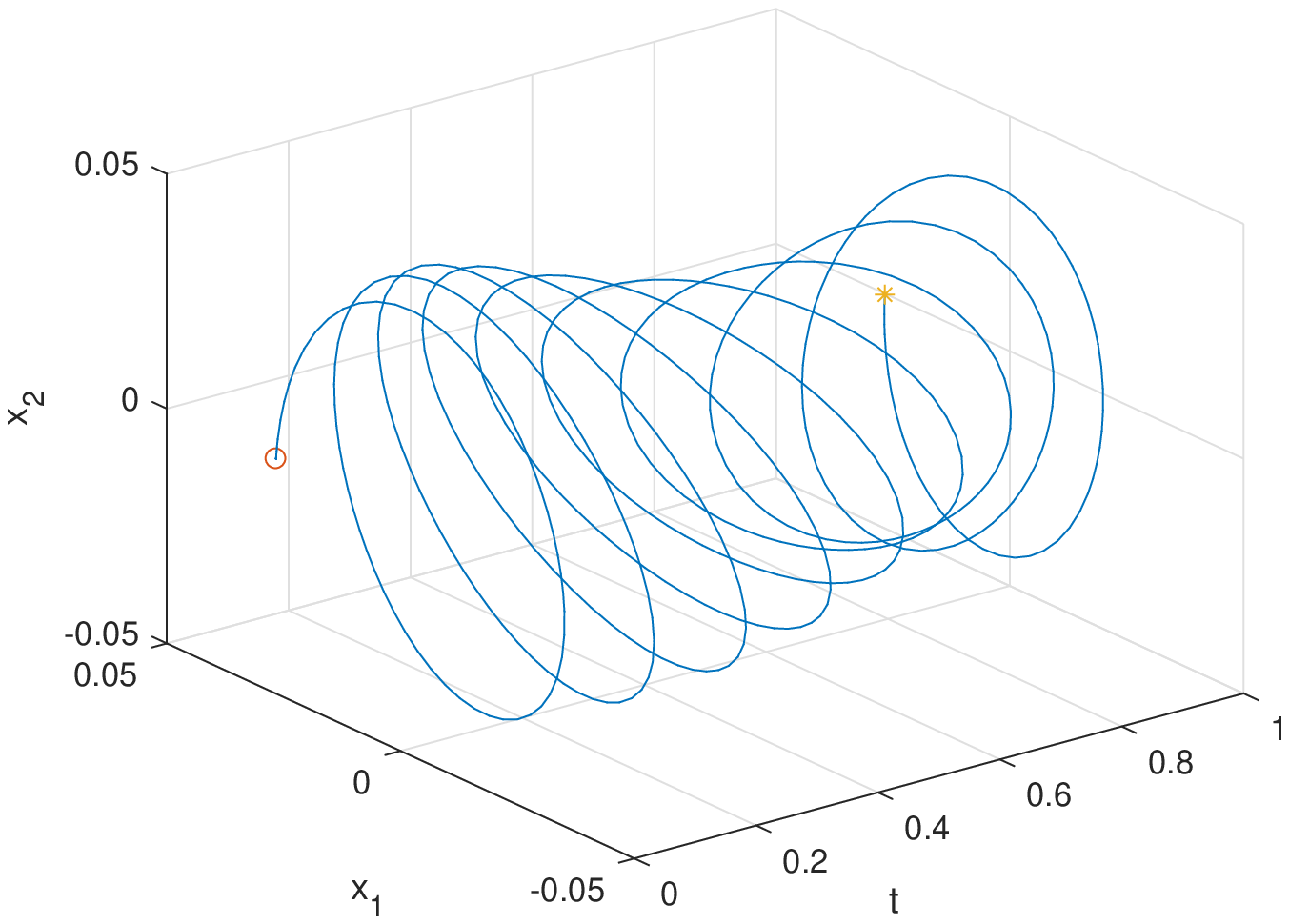}} \hspace{1cm}\resizebox{7cm}{6cm}{\includegraphics{ef3919c.eps}}
\ec
\caption{The periodic motion $(t,\bfe_{3,9/19}(t))$ and its distributions of the potential and kinetic energies. }
\lb{pon-kin3}
\end{figure}
\fi

    \bb{ex} \lb{p=5}
{\rm For $p=5$, one has
    \[
    \mcs_5\supset (2/5, 1/2).
    \]
Then $\mcs_5^*$ contains the following two numbers
    \[
    \lm = 3/7, \ (<) \ 4/9.
    \]
One has the numerical results
    \[
    \mu_{5,3/7}=0.6776 > \mu_{5,4/9}=0.5293, \andq \pi_{5,3/7}=9.3183 <\pi_{5,4/9}=12.2510
    \]
Notice that  the four ratios are
    \[
    \f{\pi_{5,3/7}}{\pi}=2.9661< 3< \f{\pi_{5,3/7}}{\pi_5}=3.3032,\q \f{\pi_{5,4/9}}{\pi}=3.8996 < 4 < \f{\pi_{5,4/9}}{\pi_5}=4.3428,
    \]
and the fifth is
    \[
    \f{\pi_{5,4/9}}{\pi_{5,3/7}}=1.3147
    \]
are all non-integers, while the five reciprocals are all $<1$.  We conclude that $\pi_{5,3/7}$  and $\pi_{5,4/9}$ can yield two new eigenvalue sequences. The new eigenvalues resulted are
    \[
    \la_{5,3/7}=(2 \pi_{5,3/7})^5 = 224819.2113, \andq \la_{5,4/9}=(2 \pi_{5,4/9})^5=883095.5120
    \]
Both of them are very large, compared with the known basic eigenvalues $(2\pi_5)^5$ and $(2\pi)^5$.\qed
}
    \end{ex}

\ifl    \bb{rem}\lb{ordering}
The ordering in \x{ord1} and \x{ord2} for index $\lm$, scaling momenta and eigenvalues are suggestive, while the ordering in \x{ra1} and \x{ra2} for ratios of eigenvalues can be explained in \x{pipr} below.
    \end{rem}
\fi

The corresponding periodic motions of Example \ref{p=5} are graphed in Figures \ref{persolp51} and \ref{persolp52}. Moreover, let us take the periodic motion $\bfe_{5,3/7}(t)$ as an example. In order that the motion is more clear, it can be compared with its potential energy and the kinetic energy
    \[
    P_{5,3/7}(t) := \la_{5,3/7} \|\bfe_{5,3/7}(t)\|^5/5, \andq K_{5,3/7}(t):= 4\|\dot \bfe_{5,3/7}(t)\|^5/5\equiv 1-P_{5,3/7}(t).
    \]
See Figure \ref{pon-kin5}. 

\begin{figure}[httb]
\bc
\resizebox{7cm}{6cm}{\includegraphics{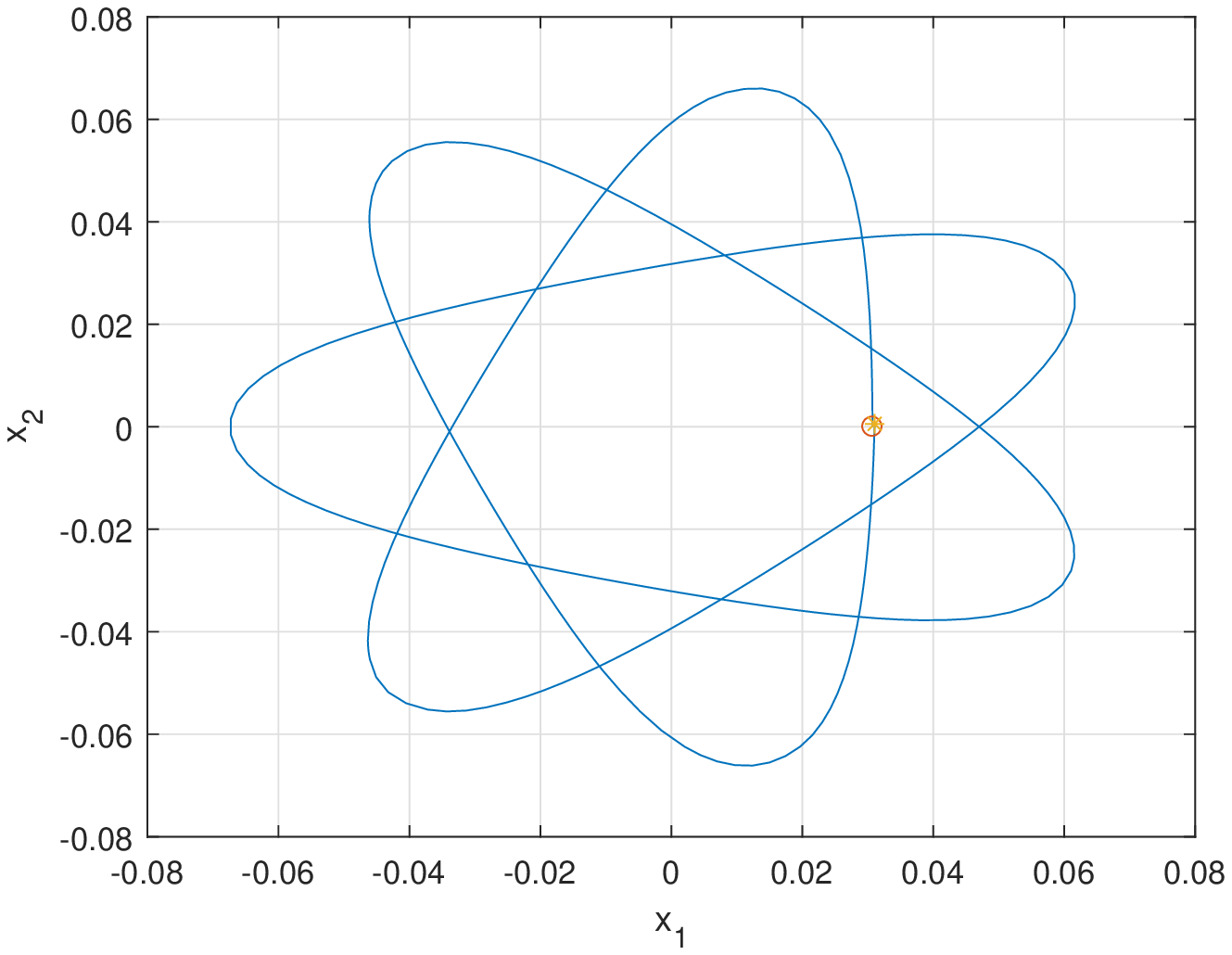}}\hspace{1cm} \resizebox{7cm}{7cm}{\includegraphics{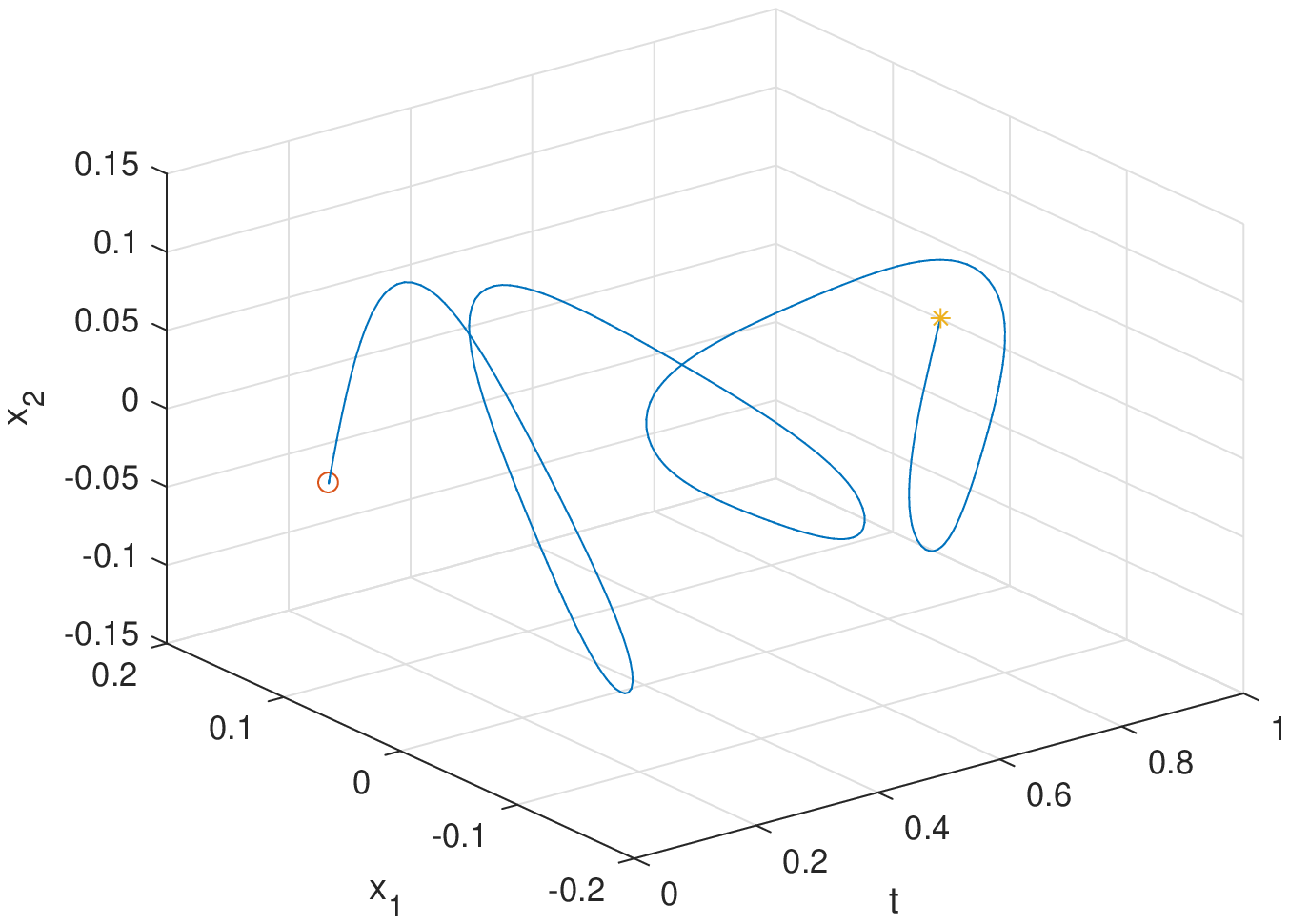}}
\ec
\caption{The periodic orbit (left) and trajectory (right) of the eigenfunction $\bfe_{5,3/7}(t)$. }
\lb{persolp51}
\end{figure}

\begin{figure}[httb]
\bc
\resizebox{7cm}{6cm}{\includegraphics{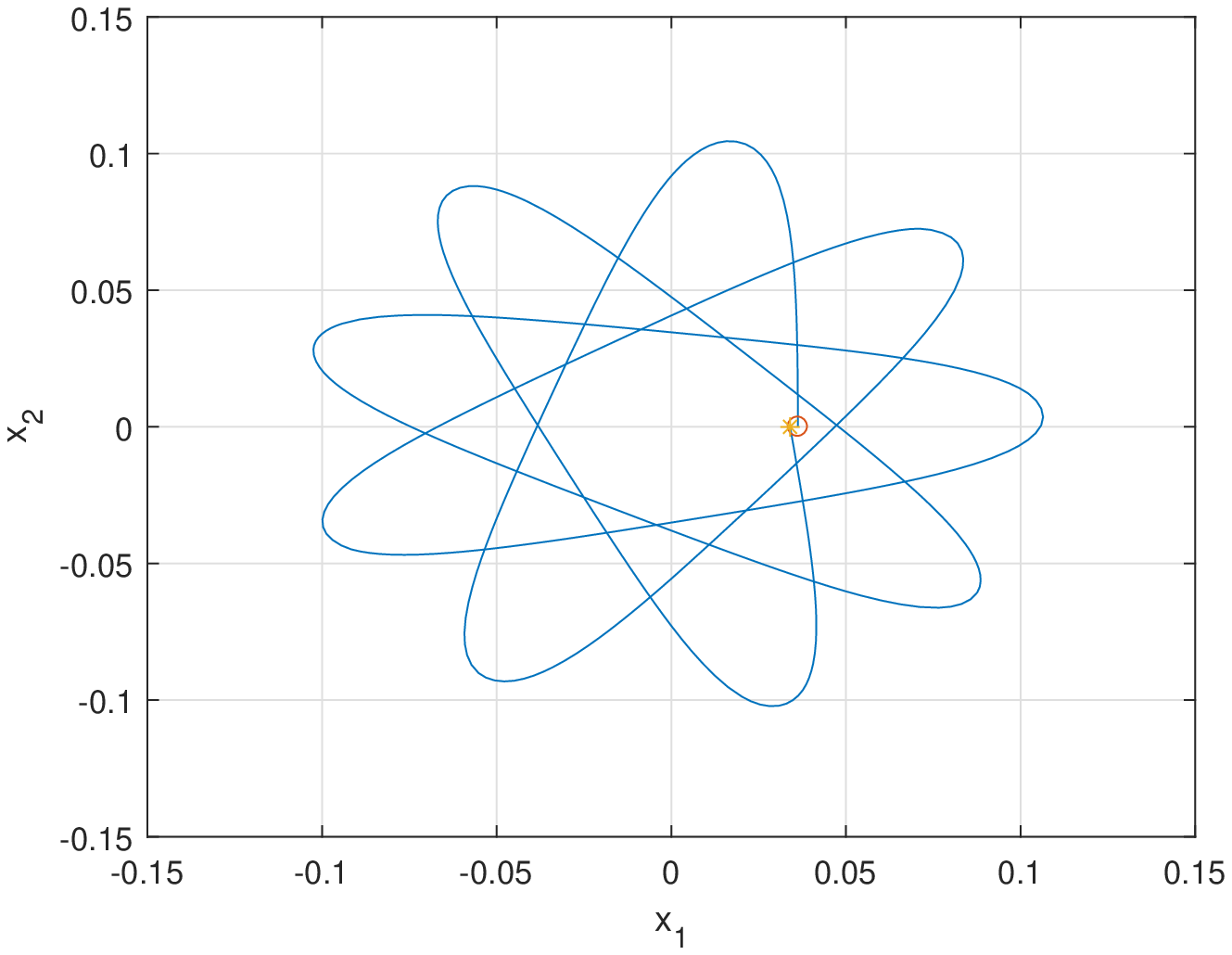}}\hspace{1cm} \resizebox{7cm}{7cm}{\includegraphics{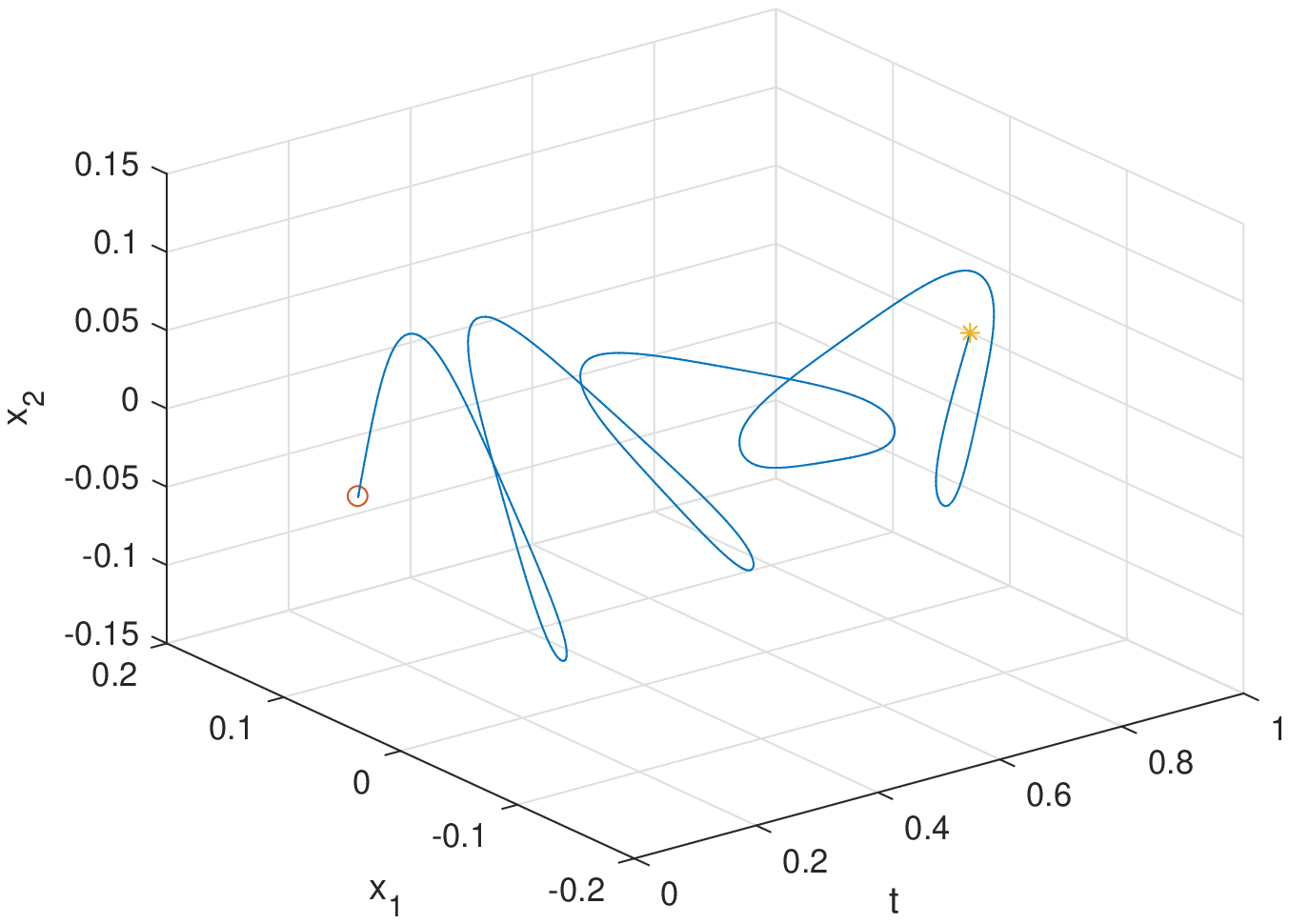}}
\ec
\caption{The periodic orbit (left) and trajectory (right) of the eigenfunction $\bfe_{5,4/9}(t)$. }
\lb{persolp52}
\end{figure}

\begin{figure}[httb]
\bc
\resizebox{7cm}{7cm}{\includegraphics{ef537b.eps}} \hspace{1cm} \resizebox{7cm}{6cm}{\includegraphics{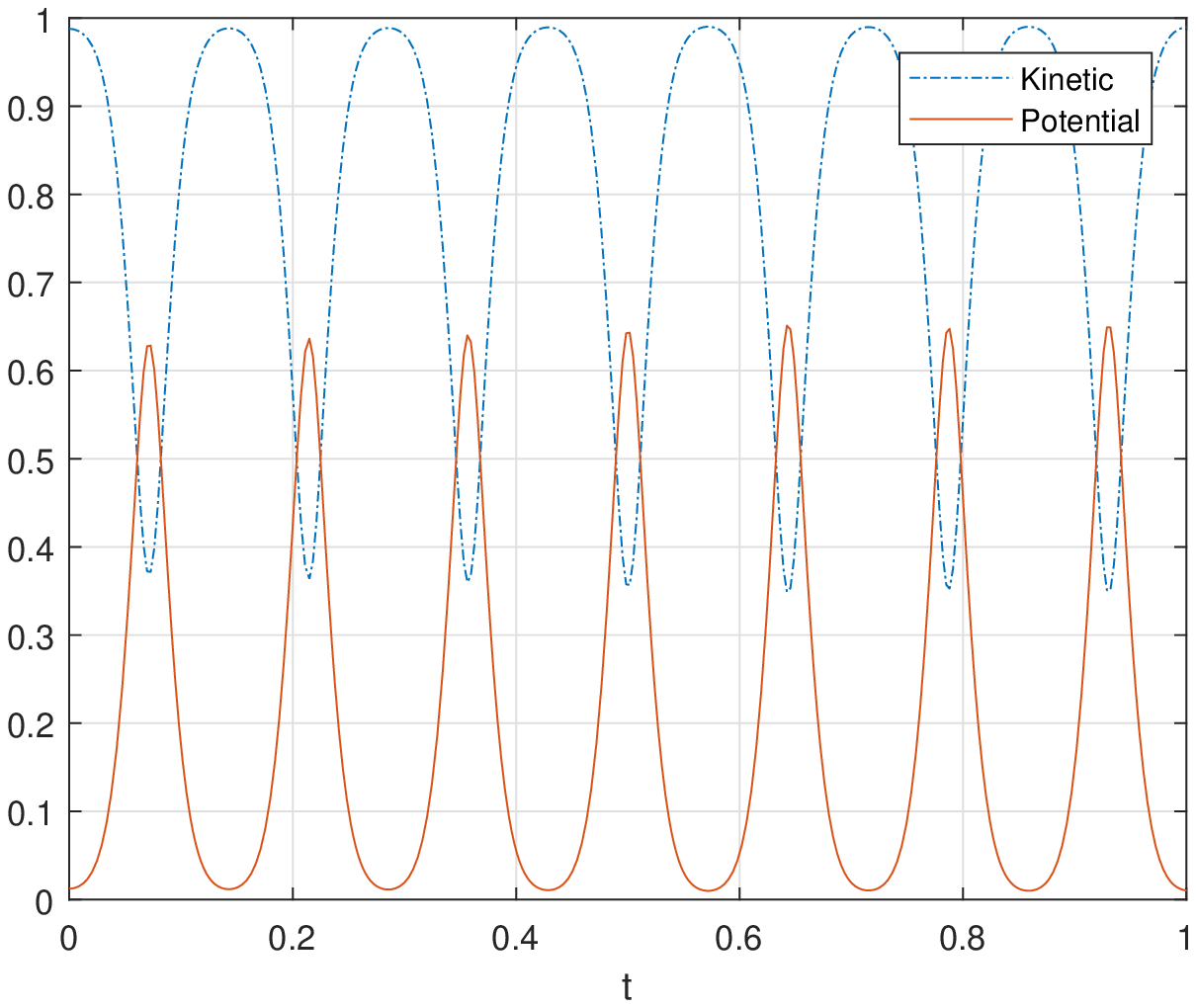}}
\ec
\caption{The periodic motion $\bfe_{5,3/7}(t)$ and its distributions of the potential and kinetic energies. }
\lb{pon-kin5}
\end{figure}

In this paper, the periodic motions $\bfe_p(t;\la_\plm)$ are indexed using $\lm\in \mcs_p^*$. From Figures \ref{persolp3}--\ref{persolp52}, the denominators $m$ and the numerators $\ell$ of $\lm$ are visible from the planar orbits and the spatial trajectories respectively. These are reasonable from the defining equalities in \x{Tdef}---\x{vpt}, and \x{Eq0} and are common in the Lissoajous figures. Moreover, the bigger $m$ and $\ell$ are, the much complicated the periodic motions $\bfe_p(t;\la_\plm)$ are. For example, one has from \x{TL} that when $p> (9+3\sqrt{5})/2= 7.8541$, $\lm=1/3\in \mcs_p^*$. Hence the $p$-Laplacian admits an eigenvalue $\la_{p,1/3}$ having a relatively simple periodic eigenfunction $\bfe_p(t; \la_{p,1/3})$.

\section{Conjectures and Discussions}
\setcounter{equation}{0} \lb{sec6}

\subsection{Conjectures}

Note that all of our constructions for eigenvalues are based on the scaling momenta of eigenfunctions. Towards to the panorama of the structure of the spectral set of the $p$-Laplacian, we impose two conjectures. Both are related with scaling momenta.

This first one is suggested by Corollary \ref{pana1}.

    \bb{conj} \lb{conj1}
Let $p\ne 2$. Problem \x{pd}-\x{bc} has infinitely many different eigenvalue sequences with the scaling momenta in any given subinterval of $(0,1)$.
    \end{conj}

This implies that the spectral set $\Si_p^*$ can admit eigenfunctions of scaling momenta which are densely distributed on $[0,1]$.

The second one is suggested by Theorem \ref{main1}. See \x{sip2}. It is possible to choose an integer-independent sequence $\{\vpi_1,\vpi_2, \dd\}$ such that
    \[
    \Si_p^*=\z\{(2n \vpi_m)^p: m , \ n\in \N \y\}.
    \]
Let us call $(2\vpi_m)^p$, $m=1,2,\dd$ the basic periodic eigenvalues.

    \bb{conj} \lb{conj2}
Let $p\ne 2$. Associated with each basic periodic eigenvalue $(2\vpi_m)^p$, $m\in \N$, all eigenfunctions have the same scaling momentum.
    \end{conj}

This means that, for basic periodic eigenvalues, the scaling momenta are dependent only on eigenvalues themselves, not on eigenfunctions.

To obtain the solutions to these conjectures, it is an important issue to understand the monotonicity of the functions $\tpm$ and $\spm$ in $\mu$, and of $\upm=\tpm/\spm$ as well. Though the reduced system \x{S2d} is integrable, we find some difficulty because \x{S2d} is not an algebraic differential system. We hope the formulas in the appendix of \cite{Levi} are helpful.

\ifl

    \bb{rem} \lb{sc-mm}
(1) Because of Example \ref{se-lan} (3) and (4), for the case $p=2$, eigenfunctions associated with the same eigenvalue may have different scaling momenta. See \x{E2t} and \x{Sme2t}. Hence the necessity parts of Proposition \ref{prop1} are stated as that {\it some}, {\it not all}, associated eigenfunctions have scaling momenta $0$ or $1$.

(2) On the other hand, whenever $p\ne 2$, it will be conjectured that {\it all} eigenfunctions associated with the same eigenvalue must have the same scaling momenta, i.e., the scaling momentum is dependent only on $\la$, not on the choice of the associated eigenfunctions. However, we are not able to prove this in the present paper. Due to this very delicate problem, the future analysis and construction of eigenvalues of the $p$-Laplacian are a little bit complicated.
    \end{rem}

Note that for the case $p=2$, we have known from Example \ref{se-lan} that scaling momenta are depending on eigenfunctions, not only on eigenvalues. This may be one of the important differences between the Laplacian and the $p$-Laplacian.

\note{Unfortunately, this fails even when $p\ne 2$. For example, when $\pi_p$ is rational, one has $\la_0:=\la_{n_0} = \tl\la_{m_0}$ for some $n_0, \ m_0\in \N$. For such an eigenvalue $\la_0$ of the $p$-Laplacian, it does admit some eigenfunctions of scaling momenta being $0$, and some other eigenfunctions of scaling momenta being $1$. Hence Conjecture \ref{conj2} may only be true for a set of exponents $p\in(1,\oo)$ of the full Lebesgue measure. In this case, one needs to exclude the rational-independence of the basic $\pi_{p,*}$'s in defining the eigenvalue sequences. Of course this is a much delicate problem. 2021.3.27.}
\fi

\subsection{Discussions} 

When $p=2$, system \x{pd} is a Newtonian equation which is just an uncoupled system of linear oscillations. The structure of the spectral set $\Si_2^*$ is simple.

When $p\ne 2$, system \x{pd} is an Euler-Lagrange equation in the Lagrangian mechanics. The cases $p<2$ and $p>2$ are usually used to model problems with spatial singularity or spatial degeneracy, respectively. The result we have proved in this paper shows that the spectral set $\Si_p^*$ is more complicated. However, from our construction and Conjecture \ref{conj2}, it seems that structure of the associated scaling momenta for basic periodic eigenvalues is simple. That is, from the point of view of scaling momenta, the spectral problems of the $p$-Laplacians $(p\ne 2)$ themselves are not degenerate in a certain sense. Moreover, as $p\to 2$, it can be expected that each branch of eigenvalues will shrink into a single eigenvalue of the Laplacian, meanwhile the scaling momenta of different branches of eigenvalues will accumulate into the whole interval $[0,1]$.

Some differences between the spectral problems of the Laplacian and the $p$-Laplacians are listed in Table \ref{tab1}.

\begin{table}
\centering
\caption{Differences between the Laplacian and the $p$-Laplacian with $d=2$} \lb{tab1}
\begin{tabular}{llll}
No. & Description & $p=2$ & $p\ne 2$ \\
1. & Equation and system & Newton Equation & Euler-Lagrange Equation \\
& & Linear, uncoupled & Nonlinear, coupled \\
2. & Modelling & Regular &  Singular for $p<2$, and \\
& &  &  degenerate for $p>2$\\
3. & Spectral set & A single sequence & Infinitely many sequences \\
4. & Scaling momenta of & $[0,1]$ & A single value $\mu=M(\la)$ \\
& a basic eigenvalue $\la$ & & (Conjectured)
\end{tabular}
\end{table}


\subsection{Comments} 

The spectral problem on different kinds of the $p$-Laplacians is a fascinating problem. Let us just mention a quite few works.

For the one-dimensional $p$-Laplacian with periodic potentials, Zhang \cite{Zh01} had used the rotation number from dynamical systems to study the periodic eigenvalues and imposed some problems.  Later, Binding and Rynne \cite{BR07} proved that, with some special choice of periodic potentials, the structure of the periodic spectral set is different from the Hill's equations \cite{MW66}.

For the $p$-Laplacian in higher dimensional domains, several interesting new phenomena have been found by Dr\'abek and his collaborators \cite{BDG, BD17}.



\section*{Acknowledgments}

The authors are partially supported by the National Natural Science Foundation of China (Grant No. 11771315 and No. 11790273).

The second author would like to express his sincere thanks to many professors for their professional helps during a long time of the preparation of this paper. The open problem was mentioned to him  many years ago by Ra\'ul Man\'asevich, one of the imposers of the problem. Jaume Llibre and Rafael Ortega offered him many helps for clarifying the complete integrability of the $p$-Laplacian and its reduction to planar dynamical systems. Jean Mawhin and Milan Tvrd\'y also gave him several useful comments. He was also benefit from discussions with many professors from China, especially with Xiuli Cen, Jibin Li, and Xiang Zhang. All of these are gratefully acknowledged.

\vfill\hfill \fbox{\small Ver. 1, 2021-4-9}

\end{document}